\documentclass[a4paper,11pt]{amsart}
\usepackage{amssymb,amsmath,amsfonts,amscd,euscript,accents}
\usepackage{hyperref}
\usepackage{color} 

\usepackage{enumitem}

\usepackage{tikz, tikz-3dplot, pgfplots}
\usepackage{tikz-cd}

\numberwithin{equation}{section}
\newtheorem{trm}{Theorem}
\newtheorem{thm}{Theorem}[section]
\newtheorem{prop}[thm]{Proposition}
\newtheorem{lem}[thm]{Lemma}
\newtheorem{cor}[thm]{Corollary}
\newtheorem{rem}[thm]{Remark}

\newtheorem{example}[thm]{Example}
\newtheorem{dfn}[thm]{Definition}

\newcommand{\bC}{{\mathbb C}}
\newcommand{\bR}{{\mathbb R}}
\newcommand{\bP}{{\mathbb P}}
\newcommand{\bQ}{{\mathbb Q}}
\newcommand{\bZ}{{\mathbb Z}}
\newcommand{\bN}{{\mathbb N}}
\newcommand{\eF}{{\EuScript{AF}}}
\newcommand{\msl}{\mathfrak{sl}}
\newcommand{\fh}{{\mathfrak h}}
\newcommand{\fb}{{\mathfrak b}}

\newcommand{\fn}{{\mathfrak n}}
\newcommand{\mb}{\mathbf}
\newcommand{\mr}{\mathrm}
\newcommand{\bdim}{\mb{dim}}

\begin{document}

%\begin{Frontmatter}

\title[Quiver Grassmannians and affine flag varieties]
{Generalized juggling patterns, quiver Grassmannians and affine flag varieties}

\author{Evgeny Feigin}
\address{E. Feigin:\newline
HSE University, Faculty of Mathematics, Usacheva 6, Moscow 119048, Russia\newline
Faculty of Mathematics and Computer Science, Weizmann Institute of
Science, POB 26, Rehovot 76100, Israel
}
\email{evgfeig@gmail.com}
\author{Martina Lanini}
\address{M. Lanini:\newline Dipartimento di Matematica\\ Universit\`a di Roma ``Tor Vergata'',  Via della Ricerca Scientifica 1, I-00133 Rome, Italy}
\email{lanini@mat.uniroma2.it}
\author{Alexander P\"utz}
\address{A. P\"utz:\newline
Faculty of Mathematics\\
Ruhr-University Bochum\\
Universit\"atsstra\ss e 150\\
44780 Bochum\\
Germany}
\email{alexander.puetz@ruhr-uni-bochum.de}


\begin{abstract}
The goal of this paper is to clarify the connection between certain structures from the theory of totally nonnegative Grassmannians, quiver Grassmannians for cyclic quivers and the theory of local models of Shimura varieties. More precisely,  we  generalize the construction from our previous paper relating the combinatorics and geometry of quiver Grassmanians to that of the totally nonnegative Grassmannians.  The varieties we are interested in serve as 
realizations of local models of Shimura varieties. We exploit quiver representation techniques to study the quiver Grassmannians of interest and, in particular, to describe explicitly embeddings into affine flag varieties which allow us to realize our quiver Grassmannians as a union of Schubert varieties therein. 
\end{abstract}

%\end{Frontmatter}
\maketitle

\iffalse
\begin{enumerate}
    \item Recall (Cyclic) Quiver Grassmannians -- Alex
    \item Properties of $X(k,n,\omega)$ -- Alex
    \item in particular: torus action, resolution of singularities -- Alex
    \item generalised juggling patterns for the fixed points in $X(k,n,\omega)$ (link to Weyl group - affine permutations- also see next item) -- Alex    
    \item geometric explanation of the property of the cells being a lower ideal (bounded affine permutations within affine permutations) - Alex
    \item Affine Flag + Grassmannian of type $\mathfrak{gl}_n$ -- Martina
    \item Schubert varieties therein -- Martina
    \item Embedding of $X(k,n,\omega)$ in the affine flag  -- Evgeny
    \item the image is the union of Schubert varieties (given explicitly)

    \item link it to Shimura varieties -- Evgeny
    \item in particular the problem with the example in our previous paper concerning flatness  -- Evgeny
    \item connection to TNN Grassmannian -- to be discussed 
    \item Are there formulae for the Poincar\'e polynomials or the KT-bases? -- to be discussed
\end{enumerate}
\fi

\section{Introduction}
Quiver Grassmannians are natural generalizations of the classical Grassmannians and flag varieties. 
In short, given a quiver $Q$, a $Q$-representation $M$ and a dimension vector ${\bf e}$ one
considers the variety ${\rm Gr}_{\bf e}(M)$ consisting of ${\bf e}$-dimensional 
subrepresentations of $M$ (\cite{CaRe2008,CI20,Scho92}). Quiver Grassmannians were extensively studied during the last two decades and proved
to be useful in various areas of mathematics (\cite{CFR12,CL,Re13}). In this paper we consider certain
quiver Grassmannians $X(k,n,\omega)$ for cyclic quivers, generalizing  
\cite{FLP21}. 
The varieties $X(k,n,\omega)$ naturally show up in the theory of local models of Shimura varieties \cite{G01,HG02,PRS13} and in the theory of totally nonnegative Grassmannians \cite{FLP21,Kn08,Lus94,Lus98a,Lus98b,Pos06}. We give some details below.

Let $\Delta_n$ be a cyclic quiver on $n$ vertices. For a positive integer 
$\omega$  we consider a $\Delta_n$ module $U_{n\omega}$ defined as follows. Let 
$A_\infty$ be an infinite (in both directions) equioriented quiver of type $A$.
Let us consider a natural $\mod n$ projection from $A_\infty$ to $\Delta_n$. Let $V(i,j)$, $i\le j$, 
be the indecomposable representation of $A_\infty$ supported on vertices from $i$ to $j$.
Then the $\Delta_n$ module $U_{n\omega}$ is obtained as the projection of the direct sum 
of $A_\infty$ modules $V(i+1,i+n\omega)$ for $i=0,\dots,n-1$. In particular, the dimension of 
$U_{n\omega}$ is equal to $(n\omega,\dots,n\omega)$.

Now let us fix $k=1,\dots,n-1$. We define $X(k,n,\omega)$ to be the $\Delta_n$  quiver Grassmannian 
${\rm Gr}_{\bf e}(U_{n\omega})$, where ${\bf e}=(k\omega,\dots,k\omega)$. The $\omega=1$ case was considered in
our previous paper \cite{FLP21}. In particular, we showed that the topological and combinatorial
properties of $X(k,n,1)$ are closely related to that of the totally nonnegative Grassmannians.
The following theorem holds.
\begin{trm}\label{trm:1}
The variety $X(k,n,\omega)$ is a projective equidimensional variety of dimension 
$\omega k(n-k)$. The number of  irreducible components is $\binom{n}{k}$. 
Each component is normal, Cohen-Macaulay, has rational singularities and 
admits a desingularization by a certain smooth quiver Grassmannian.
\end{trm}

Recall that in \cite{FLP21} a link between the varieties $X(k,n,1)$ and the totally nonnegative Grassmannians was described. More precisely, we proved that $X(k,n,1)$ admits a Bialynicki-Birula
decomposition \cite{Birula1973}, which is also a cellular decomposition. The poset of cells was
identified with the (reversed) cell poset of the corresponding tnn Grasmannian. We prove the
following generalization. 

\begin{trm}\label{trm:2}
The varieties $X(k,n,\omega)$ admit a Bialynicki-Birula decomposition with each stratum being an affine cell.
Each cell contains a unique fixed point under an appropriate algebraic torus action and the cell is the orbit of this point under the action of the automorphism group of the $\Delta_n$ module $U_{n\omega}$. The cells are labeled by a natural $\omega$-generalization of the bounded affine permutations.    
\end{trm}

The cellular decomposition is stable under the action of the above mentioned algebraic torus.
%$T$, whose fixed points are centers of the cells. We
We describe the moment graph resulting from the torus action on $X(k,n,\omega)$, and investigate the poset structure on the set of cells.  

The last part of the paper is devoted to the realization of the quiver Grassmannians $X(k,n,\omega)$
inside the affine Grassmannian of type $A$. The varieties $X(k,n,\omega)$ show up in the literature
as an explicit realization of the local models of Shimura varieties (see \cite{G01,G03,P18}). Certain properties
of these local models have been studied. In particular, the embeddings into the affine Grassmannians
were constructed. We use the techniques of quiver representations in order to prove the following.

\begin{trm}\label{trm:3}
The varieties $X(k,n,\omega)$ can be identified with the union of certain $\binom{n}{k}$-many Schubert subvarieties 
inside the type $A$ affine flag variety. The Weyl group elements corresponding to these Schubert varieties
are explicitly described. The action of the automorphism group ${\rm Aut}_{\Delta_n}(U_{n\omega})$ is identified
with the action of the Iwahori subgroup. 
\end{trm}

The link between the quiver Grassmannians and affine Schubert varieties is established via the formalism
of semi-infinite wedge spaces and Sato Grassmannians.

Finally, let us formulate a natural question which remains open. For all $\omega\ge 1$ construct
spaces ${\rm Gr}(k,n,\omega)_{\ge 0}$ which generalize  the totally nonnegative Grassmannians
${\rm Gr}(k,n)_{\ge 0}$, in the sense that ${\rm Gr}(k,n,\omega)_{\ge 0}$ admits a decomposition into cells of the form $\bR_{>0}^M$ such that the cell poset is dual
to that of $X(k,n,\omega)$. (recall that the latter is described in terms of the generalized bounded affine permutations).

The paper is organized as follows.  In Section \ref{sec:Geometric-properties} we study the geometric properties of the quiver Grassmannians $X(k,n,\omega)$; Theorem ~\ref{trm:1} is proved here. %(c.f. Theorem~\ref{thm:geometric-properties}). 
In Section \ref{sec:moment-graph+Cohomology} we describe the moment graphs of $X(k,n,\omega)$ and discuss their cohomology; in particular, we construct a cyclic group action on the  equivariant cohomology. Theorem~\ref{trm:2} is proved in Section \ref{sec:Poset-structures}. We introduce the $\omega$-generalized versions of the bounded affine permutations, juggling patterns and provide a combinatorial model for the Poincar\'e polynomials. Theorem~\ref{trm:3} is proved in Section \ref{sec:quiver-Grass-in-affine-flag} and the preliminaries on the Sato Grassmannians and affine flag varieties can be found in Section \ref{sec:Affine-flag}. In the appendix we provide a correction of a computation from \cite{FLP21}.

\section{Geometric properties}\label{sec:Geometric-properties} %of the main object? 
\subsection{Quiver Grassmannians}
In this section we recall the definition of quiver Grassmannians. For more detail on the representation theory of quivers see \cite{Schiffler2014}. A finite quiver $Q$ consists of a finite set of vertices $Q_0$, a finite set of oriented edges $Q_1$ between the vertices.  % and two maps $s,t : Q_1 \to Q_0$ providing an orientation of the edges with source $s_a$ and target $t_a$ for all $a \in Q_1$. 
A $Q$-representation $M$ is a pair of tuples, with a tuple $(M^{(i)})_{i \in Q_0}$ of $\bC$-vector spaces over the vertices, and a tuple $(M_a )_{a \in Q_1}$ containing linear maps between the vector spaces%$M^{(i)}$
, along the arrows in $Q_1$. 

A morphism $\psi$ of $Q$-representations $M$ and $N$ is a collection of linear maps $\psi_i : M^{(i)} \to N^{(i)}$ such that $\psi_j \circ M_a = N_a \circ \psi_i$ holds for all edges $a : i \to j$.  The set of all $Q$-morphisms from $M$ to $N$ is denoted by $\mathrm{Hom}_Q(M,N)$. The category of finite dimensional complex $Q$-representations is $\mathrm{rep}_\bC(Q)$. %Its dimension is abbreviated as $[M,N]:= \dim \mathrm{Hom}_Q(M,N)$.

The entries of the dimension vector $\textbf{dim} \, M \in \bZ^{Q_0}$ of a quiver representation $M \in \mathrm{rep}_\bC(Q)$ are given by $ \dim_\bC 
 M^{(i)}$ for all $i \in Q_0$. A subrepresentation $N \subseteq M$ is parameterized by a tuple of vector subspaces $N^{(i)} \subset M^{(i)}$, such that $M_a(N^{(i)}) \subseteq N^{(j)}$ holds for all arrows $a : i \to j$ of $Q$. 
\begin{dfn}
For $\mb{e} \in \bZ^{Q_0}$ and $M \in \mathrm{rep}_\bC(Q)$, the \textbf{quiver Grassmannian} $\mr{Gr}_{\mb{e}}(M)$ is the variety of all $\mb{e}$-dimensional subrepresentations of $M$.
\end{dfn}
For a point $U \in {\rm Gr}_{\bf e}(M)$ the isomorphism class $\mathcal{S}_U$ in the quiver Grassmannian is called stratum and is irreducible (cf. \cite[Lemma~2.4]{CFR12}). The automorphism group $\mr{Aut}_Q(M) \subset \mr{End}_Q(M)=\mathrm{Hom}_Q(M,M)$ acts on $\mr{Gr}_{\mb{e}}(M)$ as
\[ A.\big(U^{(i)}\big)_{i \in Q_0} := \Big( A_i\big(U^{(i)}\big) \Big)_{i \in Q_0} \quad \mr{for} \  A \in \mr{Aut}_Q(M) \ \mr{and} \ U \in {\rm Gr}_{\bf e}(M).\]
%for $A \in \mr{Aut}_Q(M)$ and $U \in {\rm Gr}_{\bf e}(M)$.

\subsection{Cyclic quiver Grassmannians}
%\begin{example}
The equioriented cycle $\Delta_n$ is the quiver with vertex set $\bZ_n := \bZ/n\bZ$ and arrows $a: i \to i+1$  for all $i \in \bZ_n$. For a $\Delta_n$-representation we write $M_i$ instead of $M_a$ for the map along the arrow $a: i \to i+1$. Now we define the $\Delta_n$-representation $U_{m}$, for $m\geq 2$: Take the vector spaces $M^{(i)}:=\bC^m$ for all $i \in \bZ_n$ and let $B^{(i)}:=\{v_j^{(i)} \, : \, j \in [m] \}$ be the standard basis of the $i$-th copy of $\bC^{m}$. Then each map $M_i$ sends $v_j^{(i)}$ to $v_{j+1}^{(i+1)}$ for $j \in [m-1]$ and $v_m^{(i)}$ to zero.  
%and linear maps (in the standard basis of $\bC^n$) of the form
%\[ M_i = \begin{pmatrix} 
%0 & 0 & \dots & 0 & 0 & 0\\
% 1 & 0 & \dots & 0 & 0 & 0\\
%  & 1 & \ddots & \vdots & \vdots & \vdots\\ 
%  & & \ddots  & 0 & 0 & 0\\
%  & 0 &  & 1 & 0 & 0\\
%  &  &  &  & 1 & 0\\
%\end{pmatrix} \] for all $a: i \to i+1$ is denoted by $U_{[n]}$.

%If we replace replace the $1$'s in the matrix above by $r\times r$ identity matrices and the zeros by $r\times r$ zero matrices we obtain a $\Delta_n$-representation with $rn$-dimensional vector spaces, which we denote by $U_{[r],n}$. The following is the main object of this section.
%\end{example}

%Let $\Delta_n$ be the cyclic quiver on $n$ vertices. 
\begin{dfn}
Let us fix numbers $k, \omega\ge 1$ with $k\leq n$. Define 
%and consider
%the $\Delta_n$ module $U_{[n],\omega}=\bigoplus_{i=1}^n U(i;\omega n)$. We consider 
the quiver Grassmannian
\[
X(k,n,\omega) := {\rm Gr}_{(k\omega,\dots,k\omega)}(U_{\omega n}).
\]
\end{dfn}
In particular, $X(k,n,1)$ coincides with the quiver Grassmannian $X(k,n)$ studied in our paper \cite{FLP21}. 

\begin{rem}\label{rem:Shimura-var}
We note that the varieties $X(k,n,\omega)$ show up in the theory of local models of
Shimura varieties (see \cite{G03,HG02,PRS13,P18,PZ13,PZ22}). More precisely, $X(k,n,\omega)$ appear as concrete realizations of the local models of Shimura varieties for $G=GL_n$ and minuscule coweights, as discussed, for example, in \cite[\S7.1]{PRS13}.  
\end{rem}

\begin{rem}
One may vary the representation $U_{\omega n}$ keeping its dimension unchanged. Then one gets a family of quiver Grassmannians in the spirit of \cite{CFFFR2017,F12}. It would be interesting to study this family.
\end{rem}

\subsection{Torus Actions}%---------------------------------------------------
The torus $\bC^*$ acts on the vector spaces of $U_{\omega n}$ with the weights $\mr{wt}(v_j^{(i)}):=j$ for all $i \in \bZ_n$ and $j \in [\omega n]$. This action extends to $X(k,n,\omega)$ by \cite[Lemma~1.1]{Cerulli2011}. 

The above $\bC^*$ action coincides with a cocharacter of an $n+1$-dimensional algebraic torus $T:=(\bC^*)^{n+1}$ which acts on the vector spaces of $U_{\omega n}$ via 
\[ \gamma.v_j^{(i)} =  \gamma_0^{j-1} \gamma_{i-j+1}v_j^{(i)}  \quad \mr{for} \ \big(\gamma_0,(\gamma_i)_{i\in \bZ_n}\big) \in T.\] 
\begin{rem}
This coincides with the torus action on $\Delta_n$-representations as defined in \cite{LaPu2020}. Hence it extends to $X(k,n,\omega)$ by \cite[Lemma~5.12]{LaPu2020}.
\end{rem}
\begin{lem}
The fixed points of the $\bC^*$ action and $T$ action on $X(k,n,\omega)$ coincide, and the number of fixed points is finite.
\end{lem}
\begin{proof}
The first part is a special case of \cite[Theorem~5.14]{LaPu2020} and the second part follows from \cite[Theorem~1]{Cerulli2011}.
\end{proof}
Now we want to introduce an explicit parametrization of the $T$-fixed points of $X(k,n,\omega)$. For $k\leq n$ we denote by $\binom{[n]}{k}$ the set of all  $k$-element subsets of $[n]$. The following definition generalizes the standard
definition of juggling patterns (see \cite{Lam16} or subsection~\ref{sec:affine-permutations}).
\begin{dfn}
For $k,n,\omega \in \bN$ with $k \leq n$, the set of $(k,n,\omega)$ juggling patterns is
\[ {\mathcal Jug}(k,n,\omega):= \Big\{ (J_i)_{i \in \bZ_n} \in \prod_{i  \in \bZ_n} \binom{[\omega n]}{ k \omega} : \tau_1(J_i\setminus\{ \omega n\})\subset J_{i+1} \  \mr{for} \ \mr{all} \ i \in \bZ_n \Big\}.\]   
\end{dfn}

\begin{lem}\label{lem:T-fixed-point-parametrisation}
The fixed points $X(k,n,\omega)^T$ are in bijection with ${\mathcal Jug}(k,n,\omega)$.
\end{lem}
\begin{proof}
It follows from \cite[Theorem~1]{Cerulli2011} that the vector spaces parameterizing the fixed points of $X(k,n,\omega)$ are spanned by subsets $P^{(i)} \subset B^{(i)}$ for $i \in \bZ_n$, where each subset has cardinality $k \omega$. This is encoded with the index sets from $\binom{[\omega n]}{ k \omega}$. Hence the condition 
\[M_i \Big( \big\langle v :  v \in P^{(i)} \big\rangle \Big) \subset \big\langle w :  w \in P^{(i+1)} \big\rangle \] 
translates to $\tau_1(J_i\setminus\{ \omega n\})\subset J_{i+1}$ where $P^{(i)} = \{ v^{(i)}_j : j \in  J_i \}$.
\end{proof}
\subsection{Geometric Properties}%-------------------------------------------------
In order to apply the desired results concerning the geometry of quiver Grassmannians for $\Delta_n$-representations, we need an alternative realization of $U_m$.
By $A_\infty$ we denote the infinite equioriented quiver of type $A$.
Let $V(i,j)$ be the indecomposable $A_\infty$-representation with vector spaces  $V(i,j)^{(k)}=\mathbb{C}$ for any $k\in [i,j]$ and maps $V(i,j)_{k}=\textrm{id}_{\mathbb{C}}$ for any $k\in[i,j-1]$. All other maps and vector spaces are zero. % and $V(i,j)_{k}=0$ otherwise. 
Let $F: A_\infty \to \Delta_n$ send $k$ to $k \mod n$, and $(a:k \to k+1)$ to $(\overline{a}: k \mod n \to k+1 \mod n)$. This induces the $\Delta_n$-representation $U_i(\ell)$%$:=F_{i\ast}(V)$
, with vector spaces $U_i(\ell)^{(j)}:=\bigoplus_{k\in F^{-1}(j)}V(i+1,i+\ell)^{(k)}$ for any $j\in \bZ_n$ and obvious linear maps. %A $\Delta_n$-representation $M$ is called nilpotent, if there exists a $N \in \bN$ such that \( M_{i+N} \circ \dots \circ M_{i+1} \circ M_i = 0 \fa i \in \bZ_n.\)

\begin{prop}
There is an isomorphism of $\Delta_n$-representations: %$U_m$ is isomorphic to the nilpotent $\Delta_n$-representation 
\[ U_m \cong \bigoplus_{i \in \bZ_n} U_i(m). \]
\end{prop}
\begin{proof}
All vector spaces of both representations have dimension $m$. Take the base change, such that the basis vector of $\bC$ in the $j$-th position of the summand $U_i(m)$ is mapped to the $i+j-1$-th basis vector of $\bC^{m}$ over the $i+j-1$-th vertex of $\Delta_n$. This induces the desired isomorphism.
\end{proof}
\begin{prop} \label{Aut}
The elements of the automorphism group ${\rm Aut}_{\Delta_n}(U_m)$
%of $U_{[n]}=\bigoplus_{i \in \mathbb{Z}_n} U(i;n)$
 are exactly the matrix tuples $ A = (A_i)_{i \in \mathbb{Z}_n}$ with
\[ 
A_i = \begin{pmatrix}
a^{(i)}_{1,1} & & & &\\
a^{(i)}_{2,1}& a^{(i-1)}_{1,1} & & &\\
\vdots & \vdots & \ddots& &\\
a^{(i)}_{m-1,1} & a^{(i-1)}_{m-2,1} & \hdots & a^{(i-m+2)}_{1,1} &\\
a^{(i)}_{m,1} & a^{(i-1)}_{m-1,1} & \hdots & a^{(i-m-2)}_{2,1} & a^{(i-m+1)}_{1,1}
\end{pmatrix} 
\]
where $a^{(i)}_{k,1} \in \mathbb{C}$ for all $i \in \mathbb{Z}_n$, $k \in [2,m]$ and $a^{(i)}_{1,1} \in \mathbb{C}^*$ for all $i \in \mathbb{Z}_n$. In particular, $\dim_\bC {\rm Aut}_{\Delta_n}(U_m)=mn$.
\end{prop}
\begin{proof}
We begin by computing the endomorphism algebra ${\rm End}_{\Delta_n}(U_m)$. By definition we have that $(E_i)_{i\in\bZ_n} \in {\rm End}_{\Delta_n}(U_m)$ if and only if
\[ E_{i+1} \tau_1 = \tau_1 E_i \quad \mathrm{for} \ \mathrm{all} \ i \in \mathbb{Z}_n. \]
This is the same as 
\begin{equation}\label{eqn:AutGp}
    E_{i+1} \tau_1(v_l^{(i)}) = \tau_1 E_i(v_l^{(i)}) \quad \mathrm{for} \ \mathrm{all} \ i \in \mathbb{Z}_n, l \in [m].
\end{equation}  
We write $a^{(i)}_{k,l}:= (E_i)_{k,l}$ for the matrix entries, so that \[E_i(v_l^{(i)})=\sum_{k=1}^m a_{k,l}^{(i)}v_k^{(i)}.\]
Then the equations \eqref{eqn:AutGp} are equivalent to 
\[
a^{(i)}_{k,l} = a^{(i+1)}_{k+1,l+1},\quad a^{(i)}_{k,n}= 0, \quad a^{(i)}_{n,l}= 0, \qquad \mathrm{for} \ \mathrm{all } \ k,l \in [m-1].
\]
From the above equations, it follows by induction on $m-l$ that $e_{k,l}^{(i)}=0$ for any $l>k$, and by induction on $l$ that  $e_{k,l}^{(i)}=e_{k+1,l+1}^{(i+1)}$.
%$e^{(i)}_{k,l} = e^{(i+1)}_{k+1,l+1}$, $e^{(i)}_{k,n}= 0$ and $e^{(i)}_{1,l+1}= 0$ for $k,l \in [n-1]$  where we write $e^{(i)}_{k,l}:= (E_i)_{k,l}$. 
 This implies that the $E_i$'s are of the claimed lower triangle form. Now we obtain the automorphism group by imposing the additional condition that all diagonal entries are invertible.
\end{proof}

The $\bC^*$ action on $X:=X(k,n,\omega)$ induces the decomposition:
\begin{equation*}%\label{eqn:BBdecomposition}
X=\bigcup_{p \in X^{\bC^*}} W_p,  \quad \quad  \mathrm{with} \quad \quad  W_p := \left\{ x \in X : \lim_{z \to 0} z.x =p \right\}.
%X(k,n,\omega)=\bigcup_{p \in X(k,n,\omega)^{\bC^*}} W_p,  \ \  \mathrm{with} \ \ W_p := \left\{ x \in X : \lim_{z \to 0} z.x =p \right\}.
\end{equation*}
%where $W_i$ is called attractive locus of $x_i$. 
We call this a BB-decomposition since decompositions of this type were first studied by Bialynicki-Birula in \cite{Birula1973}.
\begin{thm}\label{thm:geometric-properties}
For $\omega\ge 1$ and $k \leq n$ the variety $X(k,n,\omega)$ satisfies the following:
\begin{enumerate}[label=(\roman*)]
    \item it is a projective variety of dimension $\omega k(n-k)$;
    \item its irreducible components are equidimensional;
    \item the BB-decomposition is a cellular decomposition;
    \item the irreducible components are normal, Cohen-Macaulay and have rational singularities;
    \item the irreducible components $X_I(k,n,\omega)$ are labeled by the $k$-element subsets $I\subset [n]$. 
\end{enumerate}
\end{thm}
\begin{proof}
%The proof of this is analogous to the proof of \cite[Theorem~3.8]{Pue2020}: 
Part (i) is a special case of \cite[Lemma~4.9]{Pue2020}. Part (ii) and the labeling of the irreducible components as in part (v), are obtained from \cite[Lemma~4.10]{Pue2020}. The defining weight function of the $\bC^*$ action introduced above coincides with the one in \cite[Section~4.4]{Pue2020}. Hence \cite[Theorem~4.13]{Pue2020} implies that the parts of the $BB$-decomposition are cells and thus item (iii) holds. \cite[Lemma~4.12]{Pue2020} implies item (iv).
\end{proof}
\begin{rem}
Using the methods from \cite{Pue2020}, it is also possible to study ${\rm Gr}_{(q,\dots,q)}(U_{m})$, where the $\Delta_n$-representation $U_m$ is defined as above. But in this setting it is only possible to prove part $(iii)$ and $(iv)$ of the above theorem. For part $(i)$, $(ii)$ and $(v)$ it is crucial that $n$ divides $m$ and that $m/n$ divides $q$. Otherwise, the irreducible components are not of the same dimension, their parametrization is unknown and hence there is no dimension formula. %TODO: give example!
\end{rem}
%Here is an explicit realization of the representation $U_{[n],l}$. We fix bases $v_j^{(i)}$  of the spaces $U^{(i)}_{[n],l}$  space We fix vectors $v_1,\dots,v_{ln}$ spanning a vector space $E$  and we denote by $\tau_1: E\to E$ the linear map sending $e_i$ to $e_{i+1}$ (and $\tau_1e_{ln}=0$). We identify each space $U^{(i)}_{[n],l}$ with $E$ and all the maps with $\tau_1$.
%\martina{The following Thm refers to some torus action, but it seems to me that we have not yet defined it.}
\begin{thm}\label{trm:cells-are-strata}
Each cell $C \subset X(k,n,\omega)$ is $T$-stable and contains exactly one $T$-fixed point $P_C$. The $\mr{Aut}_{\Delta_n}(U_{\omega n})$-orbit and stratum of $P_C$ coincide with $C$.
\end{thm}
\begin{proof}
The first part follows from \cite[Theorem~5.7]{LaPu2020} in combination with \cite[Theorem~5.14]{LaPu2020}. The bijection between strata and $\mr{Aut}_{\Delta_n}(U_{\omega n})$-orbits follows from \cite[Lemma~2.28]{Pue2019}. It follows from the explicit description of the $T$-fixed points in Lemma~\ref{lem:T-fixed-point-parametrisation} that they are pairwise not isomorphic. Hence they belong to different strata.
\end{proof}
\begin{rem}
Observe that in general the $\mr{Aut}_Q(M)$-orbits in a quiver Grassmannian $\mr{Gr}_{\mb{e}}(M)$ are not cells and contain more than one torus fixed point. This is very special to the $X(k,n,\omega)$'s. Otherwise it already fails in small examples:

Consider the $\Delta_2$-representation $M$ with $M^{(1)}=M^{(2)}=\bC^2$, $M_1 = \mr{id}$ and $M_2=0$. For $e=(1,1)$ let $\bC^*$ act on the quiver Grassmannian $\mr{Gr}_{\mb{e}}(M)$ induced by the weight function $\mr{wt}(e_i)=i$, where $\{e_1,e_2\}$ is the standard basis of $\bC^2$. It follows from \cite[Theorem~1]{Cerulli2011} that $\mr{Gr}_{\mb{e}}(M)$ has two isomorphic $\bC^*$-fixed points. Hence they live in the same $\mr{Aut}_{\Delta_2}(M)$-orbit.
\end{rem}
\begin{rem}\label{rem:cell-closures}
Theorem~\ref{trm:cells-are-strata} implies that every cell closure in $X(k,n,\omega)$ is the union of smaller cells. The bijection between cells and generalized juggling patterns allows to make this description explicit (Corollary~\ref{cor:about-cells-and-their-closures}). But first we have to describe the moment graph for the $T$ action on $X(k,n,\omega)$ (Lemma~\ref{lem:moment-graph-edges}).
\end{rem}

\subsection{Desingularization}%-----------------------------
In this subsection we give a desingularization of $X(k,n,\omega)$ following the general construction in \cite{PuRe2022} (see also \cite{CFR13,FF13,KS14,S17}).

For $I \in \binom{[n]}{k}$ define the $\Delta_n$-representation
\[ U_I := \bigoplus_{i \in I}U_i(\omega n). \]
\begin{prop}
For all $I \in \binom{[n]}{k}$, the closed stratum $\overline{\mathcal{S}_{U_I}}$ is an irreducible component of
$X(k,n,\omega)$ and all irreducible components are of this form.
\end{prop}
\begin{proof}
The parametrization for the representatives of the top dimensional strata is obtained in the proof of \cite[Lemma~4.10]{Pue2020}. The second part follows from Theorem~\ref{thm:geometric-properties}.(v).
\end{proof}

Let $\hat{\Delta}_n$ be the quiver with vertex set
\[ \big\{ (i,k) \ : \ i \in \bZ_n \ \mr{and} \ k \in [\omega n] \big\} \]
and arrows
\begin{align*}
&\big\{ a_{i,k} : (i,k)\to (i,k+1) \ : \ i \in \bZ_n \ \mr{and} \ k \in [\omega n-1] \big\} \cup\\
&\big\{ b_{i,k} : (i,k)\to (i+1,k-1) \ : \ i \in \bZ_n \ \mr{and} \ k \in [\omega n]\setminus \{1\} \big\}.
\end{align*}
%We define the map $ \Lambda : \mr{rep}_\bC(\Delta_n,\mr{I}_N) \to \mr{rep}_\bC(\hat{\Delta}_n,\hat{\mr{I}}_{n,N})$ by
For $M \in \mr{rep}_\bC({\Delta}_n)$ define the $\hat{\Delta}_n$-representation %$\hat{U}_{\omega n}$
 \[ \hat{M} := \big( (\hat{M}^{(i,k)})_{i \in \bZ_n, k \in [N] }, (\hat{M}_{a_{i,k}}, \hat{M}_{b_{i,k+1}})_{i \in \bZ_n, k \in [N-1]} \big)\]
 with 
 \begin{align*}
 \hat{M}^{(i,1)} &:= M^{(i)}   &\mr{for} \ k =1\\
\hat{M}^{(i,k)} &:= M_{i+k-2} \circ  M_{i+k-3} \circ  \dots \circ M_{i+1} \circ  M_{i} (M^{(i)})  &\mr{for} \ k \geq 2\\
\hat{M}_{a_{i,k}} &:= M_{i+k-1}  &\mr{for} \ k \geq 1\\
\hat{M}_{b_{i,k}} &:= \iota : \hat{M}^{(i,k)}  \hookrightarrow \hat{M}^{(i+1,k-1)}  &\mr{for} \ k \geq 2
 \end{align*}
Here the inclusion maps along $b_{i,k}$ arise naturally from the definition of the vector spaces of $\hat{M}$.

Every $W \in \mr{rep}_\bC(\hat{\Delta}_n)$ restricts to a $\Delta_n$-representation 
\[ \mr{res}  \, W := \Big( \big(W^{(i,1)}\big)_{i \in \bZ_n}, \big(W_{b_{i,2}} \circ W_{a_{i,1}}\big)_{i \in \bZ_n} \Big).\] For $I \in \binom{[n]}{k}$ set $\hat{X}_I(k,n,\omega):= \mr{Gr}_{\bdim \, \hat{U}_I} \big(\hat{U}_{\omega n} \big)$ define the map
\[ \pi_I : \hat{X}_I(k,n,\omega) \longrightarrow X(k,n,\omega) \]
by $\pi_I(V) := \mr{res} V$ for all $V \in \hat{X}_I(k,n,\omega)$.
\begin{rem}
The vector spaces of $\hat{U}_{\omega n}$ are spanned by subsets of the bases for the vector spaces of $U_{\omega n}$. Hence the $T$ action on $U_{\omega n}$ extends to the quiver Grassmannians $\hat{X}_I(k,n,\omega)$ for $I \in \binom{[n]}{k}$ in the obvious way. The same holds for the $\bC^*$ action.
\end{rem}
The following result is a special case of \cite[Theorem~3.18, Lemma~5.3]{PuRe2022}. %TODO Put the correct reference to PR23! 
\begin{thm}\label{thm:desing}
The map
\[ \pi := \bigsqcup_{I \in \binom{[n]}{k}} \pi_I  \ : \ \bigsqcup_{I \in \binom{[n]}{k}} \hat{X}_I(k,n,\omega) \ \longrightarrow \  X(k,n,\omega)\]
 is a $T$-equivariant desingularization of $X(k,n,\omega)$.
\end{thm}

\begin{thm}\label{trm:tower-of-grassmann-bundles}
For each $I \in \binom{[n]}{k}$ the quiver Grassmannian $\hat{X}_I(k,n,\omega)$ is isomorphic to a tower of fibrations 
\[  \hat{X}_I(k,n,\omega) = X_1 \to X_{2} \to \dots \to X_{\omega n } = \mr{pt} \]%\prod_{i \in \bZ_n}  \mr{Gr}_{\hat{n}^{(i,\omega n)}}\Big(\bC^{\hat{m}^{(i,\omega n)}}\Big) \]
where %$\hat{\mb{n}} := \bdim \, \hat{U}_I$ and $\hat{\mb{m}} := \bdim \, \hat{U}_{\omega n}$ and 
each map $X_k \to X_{k+1}$ for $k \in [\omega n-1]$ is a fibration with fiber isomorphic to a product of Grassmannians of subspaces.
\end{thm}
This result is a special case of \cite[Theorem~3.21]{PuRe2022} and generalizes \cite[Theorem~7.10]{FLP21}. %TODO Put the correct reference to PR22! 
%-----------------------------
\subsection{Properties of the Desingularization}
%-----------------------------
\begin{lem}(c.f. \cite[Theorem~5.5]{PuRe2022})\label{lem:cell-decomp-desing}
For $I \in \binom{[n]}{k}$ the $T$-fixed points of $\hat{X}_I(k,n,\omega)$ are exactly the preimages of the $T$-fixed points of $X_I(k,n,\omega) \subset X(k,n,\omega)$ under $\pi_I$ (where $X_I(k,n,\omega):=\overline{\mathcal{S}_{U_I}}$). The $\bC^*$-attracting sets of these points provide a cellular decomposition of $\hat{X}_I(k,n,\omega)$.
\end{lem}

\begin{prop}\label{prop:aut-group-desing}
The automorphism group of $ \hat{U}_{\omega n}$ satisfies
\[  {\rm Aut}_{\hat{\Delta}_n}\big( \hat{U}_{\omega n}\big) \cong {\rm Aut}_{\Delta_n}\big(U_{\omega n}\big). \]
\end{prop}

\begin{proof}
From the compositions $\beta_{i+1,2} \circ \alpha_{i,1}$ for all $i \in \mathbb{Z}_n$ we obtain the same relations on each matrix $A^{(i,1)}$ of $A \in {\rm Aut}_{\hat{\Delta}_n}( \hat{U}_{\omega n})$ as for the matrix $B^{(i)}$ of $B \in {\rm Aut}_{\Delta_n}(U_{\omega n})$ (see Proposition~\ref{Aut}). Now it follows from the construction of $\hat{U}_{\omega n}$ that all other components $A^{(i,r)}$ are the lower diagonal blocks of size $\omega n-r+1$ in the matrices $A^{(i,1)}$. This implies the desired isomorphism. 
\end{proof} 

\begin{lem}
The strata in the quiver Grassmannian $\hat{X}_I(k,n,\omega)$ are exactly the ${\rm Aut}_{\hat{\Delta}_n}\big( \hat{U}_{\omega n}\big) $ orbits of the $T$-fixed points and coincide with their $\bC^*$-attracting sets.
\end{lem}

\begin{proof}
The representation $\hat{U}_{\omega n}$ is an injective bounded $\hat{\Delta}_n$ representation. Hence we can apply \cite[Lemma~2.28]{Pue2019}, to conclude that all strata are ${\rm Aut}_{\hat{\Delta}_n}\big( \hat{U}_{\omega n}\big) $ orbits and vice versa. 
It follows from the explicit description of the $T$-fixed points of $\hat{X}_I(k,n,\omega)$ from Lemma~\ref{lem:cell-decomp-desing} that they are pairwise isomorphic. This implies that each ${\rm Aut}_{\hat{\Delta}_n}\big( \hat{U}_{\omega n}\big) $ orbit contains exactly one $T$-fixed point. Hence it has to coincide with the $\bC^*$-attracting set of that fixed point.
\end{proof} 

\section{Moment Graph and Cohomology}\label{sec:moment-graph+Cohomology}%-----------------------------
\subsection{Moment Graph}\label{sec:moment-graph}%-----------------------------
There is a combinatorial object called moment graph which captures the structure of fixed points and one-dimensional orbits for suitable torus actions on complex projective varieties. The structure of this graph helps to understand the equivariant geometry of the variety. Before we describe it for $X(k,n,\omega)$ we recall the definition and some required terminology.

Let $X$ be complex projective algebraic variety $X$ acted upon by a torus $T$ with finitely many fixed points and finitely many one-dimensional $T$ orbits (i.e. the action is skeletal). The definition below is specialized to the setting that $X$ admits a $T$-stable cellular decomposition (as in our case).
\begin{dfn}Let $T$ be an algebraic torus and let $X$ be a complex projective algebraic $T$-variety. Assume that $X$ admits a $T$-stable cellular decomposition where every cell has exactly one fixed point. If the action of $T$ on $X$ is skeletal, then the corresponding moment graph $\mathcal{G}(X,T)$ is given by 
\begin{itemize}
    \item the vertex set is the fixed point set: $\mathcal{V}=X^T$;
    \item there is an edge $x\to y$ if and only if $x$ and $y$ belong to the same one dimensional $T$ orbit closure $\overline{\mathcal{O}_{x\to y}}$ and $y$ belongs to the closure of the cell containing $x$;
    \item the label of the edge $x\to y$ is the character $\alpha\in {\rm Hom}(T,\bC^*)$ the torus acts by on $\mathcal{O}_{x\to y}$. 
\end{itemize}
\end{dfn}

The edge labels are only well defined up to a sign, but since this does not play any role in the applications (e.g. computation of equivariant cohomology), we assume the labels to be fixed once and for all, and forget about this ambiguity.

The following parametrization of the $T$-fixed points of $X(k,n,\omega)$ helps to describe the structure of the one-dimensional $T$-orbits.
\begin{prop}\label{prop:bijection_jug-cells}
For $k,n,\omega \in \bN$ with $k \leq n$, there is a bijection between ${\mathcal Jug}(k,n,\omega)$ and 
\[ {\mathcal C}_{k,n,\omega} := \Big\{ (\ell_j)_{j \in \bZ_n} \in [0,\omega n]^{\bZ_n} : \bdim \bigoplus_{j \in \bZ_n} U_j(\ell_j) = (k\omega, \dots, k\omega) \in \bN^{\bZ_n} \Big\}. \]
\end{prop}
\begin{proof}
Let $\varphi : {\mathcal Jug}(k,n,\omega) \to \mathcal{C}_{k,n,\omega}$ send $\mathcal{J}_\bullet = (J_i)_{i \in \bZ_n}$ to $\ell_\bullet = (\ell_j)_{j \in \bZ_n}$ with
\[ \ell_j := \max \big(\{ r \in [\omega n] : \omega n -r+1 \in J_{j-r+1} \} \cup \{ 0 \} \big). \]
%\begin{cases} \max \{ r \in [\omega n] : \omega n -r+1 \in J_{j-r+1} \}, & a\notin I,\\ a+n, & a\in I\end{cases}\]
It follows immediately from the definition of $U_{\ell_\bullet} := \bigoplus_{j \in \bZ_n} U_j(\ell_j)$ that $\bdim \, U_{\ell_\bullet} = (k\omega, \dots, k\omega)$ since each $J_i$ contains $k\omega$-many elements.

The inverse map $\varphi^{-1}: \mathcal{C}_{k,n,\omega} \to {\mathcal Jug}(k,n,\omega)$ sends $\ell_\bullet$ to $\mathcal{J}_\bullet$ where for each $\ell_j \neq 0$ and $s \in [\ell_j]$ the set $J_{j-s+1}$ contains the element $\omega n-s+1$. Clearly, the $\ell_j$'s contribute to each $J_i$ exactly $k\omega$-many times since $\bdim \, U_{\ell_\bullet} = (k\omega, \dots, k\omega)$.
\end{proof}
\begin{example}\label{ex:T-fixed-points}
For $k=n=\omega=2$, $X(k,n,\omega)$ has five $T=(\bC^*)^{2+1}$-fixed points labeled by the tuples 
\( (2,2),(3,1),(1,3),(4,0),\ \mathrm{and} \ (0,4). \)
\end{example}
Now we describe certain cut and paste moves on the segments of the elements in $\mathcal{C}_{k,n,\omega}$.
For every element $\ell_\bullet \in \mathcal{C}_{k,n,\omega}$ there are maps of the form $f_{i,j,r} : \mathcal{C}_{k,n,\omega} \to \mathcal{C}_{k,n,\omega}$ with
\[  \Big(f_{i,j,r}\big(\ell_\bullet\big)\Big)_s := \begin{cases} \ell_s & s \notin \{i,j\}\\ \ell_i-r & s=i\\  \ell_j+r & s=j\end{cases}.\]
whenever $r \in [0,\min\{\ell_i,\omega n-\ell_j\}]$ and $i-\ell_i = j - \ell_j-r \mod n$. It is straightforward to check that $f_{i,j,r}\big(\ell_\bullet\big)$ is again an element of $\mathcal{C}_{k,n,\omega}$. These cut and paste moves describe all one-dimensional $T$-orbits.
\begin{lem}\label{lem:moment-graph-edges}
    The vertices of the moment graph for the action of the torus $T$ on $X(k,n,\omega)$ are labelled by the elements of ${\mathcal C}_{k,n,\omega}$. There is an oriented edge in the moment graph from $\ell_\bullet$ to $f_{i,j,r}\big(\ell_\bullet\big)$ if and only if $\ell_i > \ell_j+r$. The label of the edge $\ell_\bullet \to f_{i,j,r}\big(\ell_\bullet\big)$ is $\epsilon_j-\epsilon_i +\delta\cdot(\ell_i-\ell_j-r)$, where $\delta(\gamma):=\gamma_0$ and $\epsilon_i(\gamma):=\gamma_i$ for any $i \in [n]$ and $\gamma=(\gamma_0,\gamma_1,\dots,\gamma_n) \in T$.
\end{lem}
\begin{proof}
    This is a special case of the description of the edges in the moment graph and their labels as given in \cite[Theorem~6.15]{LaPu2020}, translated to the description of the $T$-fixed points from Proposition~\ref{prop:bijection_jug-cells}.
\end{proof}
\begin{example}\label{ex:one-dim-orbits}
    The one-dimensional $T$-orbits between the $T$-fixed points from Example~\ref{ex:T-fixed-points} are captured in the following graph:
    
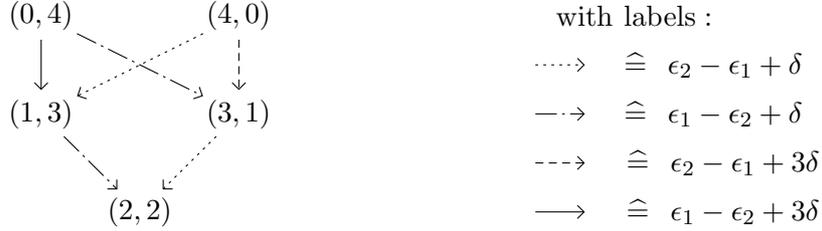
\begin{center}
        \begin{tikzpicture}[scale=1.3]
\node at (5,2) {$\mr{with}\ \mr{labels:}$};
\draw[arrows={-angle 90},dash pattern={on 1pt off 2pt on 1pt off 2pt}]  (4.0,1.5) -- (4.5,1.5);
%\draw[arrows={-angle 90}]  (4.0,1.5) -- (4.5,1.5);
\node at (5.8,1.5) {$\widehat{=} \ \ \epsilon_2-\epsilon_1+\delta$};
\draw[arrows={-angle 90},dash pattern={on 8pt off 2pt on 1pt off 2pt}]  (4.0,1.0) -- (4.5,1.0);
\node at (5.8,1.0) {$\widehat{=} \ \ \epsilon_1-\epsilon_2+\delta$};
\draw[arrows={-angle 90},dash pattern={on 3pt off 2pt on 3pt off 2pt}]  (4.0,0.5) -- (4.5,0.5);
\node at (5.9,0.5) {$\widehat{=} \ \ \epsilon_2-\epsilon_1+3\delta$};
\draw[arrows={-angle 90}]  (4.0,0.0) -- (4.5,0.0);
\node at (5.9,0.0) {$\widehat{=} \ \ \epsilon_1-\epsilon_2+3\delta$};
%\draw[arrows={-angle 90},dash pattern={on 6pt off 2pt on 2pt off 2pt}]  (4.0,-0.5) -- (4.5,-0.5); \node at (5.9,-0.5) {$\widehat{=} \ \ \epsilon_1-\epsilon_2+3\delta$};

\node at (0,0) {$(2,2)$}; 
\node at (1,1) {$(3,1)$};
\node at (-1,1) {$(1,3)$};
\node at (1,2) {$(4,0)$};
\node at (-1,2) {$(0,4)$};

\draw[arrows={-angle 90}, shorten >=9, shorten <=9,dash pattern={on 3pt off 2pt on 3pt off 2pt}]  (1,2) -- (1,1);
è\draw[arrows={-angle 90}, shorten >=15, shorten <=15,dash pattern={on 1pt off 2pt on 1pt off 2pt}]  (1,2) -- (-1,1);
\draw[arrows={-angle 90}, shorten >=15, shorten <=15,dash pattern={on 8pt off 2pt on 1pt off 2pt}]  (-1,2) -- (1,1);
\draw[arrows={-angle 90}, shorten >=9, shorten <=9]  (-1,2) -- (-1,1);
\draw[arrows={-angle 90}, shorten >=12, shorten <=12,dash pattern={on 1pt off 2pt on 1pt off 2pt}]  (1,1) -- (0,0);
\draw[arrows={-angle 90}, shorten >=12, shorten <=12,dash pattern={on 8pt off 2pt on 1pt off 2pt}]  (-1,1) -- (0,0);
            
        \end{tikzpicture}
    \end{center}
    Observe that there is no edge between $(4,0)$ and $(0,4)$ because $1 \neq 2 \mod 2$. For the same reason there can't be an edge between $(3,1)$ and $(1,3)$, and $(4,0)$ or $(0,4)$ and $(2,2)$. The label of the edge $(4,0) \to (1,3)$ is $\epsilon_2-\epsilon_1+\delta$, since $(1,3)=f_{1,2,3}(4,0)$. All other labels are computed in the same way.
\end{example}

\subsection{$T$-equivariant cohomology}\label{sec:T-equi-cohomology}%-------------------
By \cite[Theorem~6.6]{LaPu2020}, we can use the structure of the moment graph $\mathcal{G}$ as described in Lemma~\ref{lem:moment-graph-edges} to compute the ($T$-equivariant) cohomology ring of $X:=X(k,n,\omega)$. Let $R:=\mathbb{Q}[\epsilon_1, \ldots, \epsilon_n, \delta]$ and consider it as a $\mathbb{Z}$-graded ring with grading induced by $\deg(\epsilon_i)=\deg(\delta)=2$ for all $i\in [n]$. By $\alpha(\ell_\bullet,\ell'_\bullet)$ we denote the label of the edge $\ell_\bullet \to \ell'_\bullet$. \cite[Theorem~1.2.2]{GKM1998} %\cite[Theorem~1.20]{LaPu2020} 
gives the following result.
\begin{cor}\label{cor:GKM}
There is an isomorphism of ($\mathbb{Z}$-graded) rings
\[
H_T^\bullet\big(X,\mathbb{Q}\big)\simeq \left\{\big(z_{\mathcal{\ell_\bullet}}\big)_{\ell_\bullet \in\mathcal{C}({k,n,\omega})} \in\bigoplus_{\ell_\bullet \in\mathcal{C}({k,n,\omega})}R \, \Bigg\vert   \begin{array}{c}
  z_{\ell_\bullet}\equiv z_{\ell_\bullet'}\mod \alpha(\ell_\bullet,\ell'_\bullet)\\
\mathrm{for} \ \mathrm{every} \ \mathrm{edge} \ \ell_\bullet\rightarrow \ell_\bullet' \end{array}\right\}.
\]
\end{cor}
\begin{rem}
By \cite[Theorem~3.22]{LaPu2021}, $H_T^\bullet(X(k,n,\omega),\bQ)$ admits a very nice basis as a free module over $R$, namely a so-called Knutson-Tao (KT) basis (cf. \cite[Definition~3.2]{LaPu2021}). 
\end{rem}
\begin{example}\label{ex:KT-classes}
    For instance, in the case the moment graph is the one from Example \ref{ex:one-dim-orbits}, if we denote by $\alpha=\epsilon_1-\epsilon_2$, the KT classes are the following:
    
    \vspace{5mm}
    \hspace{4mm}
          % \begin{center}
        \begin{tikzpicture}[scale=0.8]
        \node at (-2,1) {$\xi_{(0,4)}=$};
\node at (0.2,0) {$0$}; 
\node at (-0.7,1) {$0$};
\node at (1.2,1) {$0$};
\node at (1.2,2) {$0$};
\node at (-1.7,2.1) {$(\alpha+\delta)(\alpha+3\delta)$};
\draw[arrows={-angle 90}, shorten >=9, shorten <=9]  (1.2,2) -- (1.2,1);
è\draw[arrows={-angle 90}, shorten >=15, shorten <=15]  (1.2,2) -- (-1,1);
\draw[arrows={-angle 90}, shorten >=15, shorten <=15]  (-1,2) -- (1.2,1);
\draw[arrows={-angle 90}, shorten >=9, shorten <=9]  (-1,2) -- (-1,1);
\draw[arrows={-angle 90}, shorten >=12, shorten <=12]  (1.2,1) -- (0.2,0);
\draw[arrows={-angle 90}, shorten >=12, shorten <=12]  (-1,1) -- (0.2,0);          
        %\end{tikzpicture}     
        %\quad
        %\begin{tikzpicture}[scale=0.8]
        \node at (5,1) {$\xi_{(4,0)}=$};
\node at (7,0) {$0$}; 
\node at (6,1) {$0$};
\node at (8,1) {$0$};
\node at (6,2) {$0$};
\node at (9,2) {$(-\alpha+\delta)(-\alpha+3\delta)$};
\draw[arrows={-angle 90}, shorten >=9, shorten <=9]  (8,2) -- (8,1);
è\draw[arrows={-angle 90}, shorten >=15, shorten <=15]  (8,2) -- (6,1);
\draw[arrows={-angle 90}, shorten >=15, shorten <=15]  (6,2) -- (8,1);
\draw[arrows={-angle 90}, shorten >=9, shorten <=9]  (6,2) -- (6,1);
\draw[arrows={-angle 90}, shorten >=12, shorten <=12]  (8,1) -- (7,0);
\draw[arrows={-angle 90}, shorten >=12, shorten <=12]  (6,1) -- (7,0);          
        \end{tikzpicture}
      %  \end{center}
      %\begin{center}

       \vspace{3mm}
    \hspace{2mm}
        \begin{tikzpicture}[scale=0.8]
        \node at (-2.9,1) {$\xi_{(1,3)}=$};
\node at (0,0) {$0$}; 
\node at (1,1) {$0$};
\node at (-1.3,1) {$\alpha+\delta$};
\node at (1,2) {$-\alpha+3\delta$};
\node at (-1.3,2) {$\alpha+\delta$};

\draw[arrows={-angle 90}, shorten >=9, shorten <=9]  (1,2) -- (1,1);
è\draw[arrows={-angle 90}, shorten >=15, shorten <=15]  (1,2) -- (-1,1);
\draw[arrows={-angle 90}, shorten >=15, shorten <=15]  (-1,2) -- (1,1);
\draw[arrows={-angle 90}, shorten >=9, shorten <=9]  (-1,2) -- (-1,1);
\draw[arrows={-angle 90}, shorten >=12, shorten <=12]  (1,1) -- (0,0);
\draw[arrows={-angle 90}, shorten >=12, shorten <=12]  (-1,1) -- (0,0);          
       % \end{tikzpicture}     
        %\quad
        %\begin{tikzpicture}[scale=0.8]
        \node at (4.7,1) {$\xi_{(3,1)}=$};
\node at (7,0) {$0$}; 
\node at (6,1) {$0$};
\node at (8,1) {$-\alpha+\delta$};
\node at (6,2) {$\alpha+3\delta$};
\node at (8,2) {$-\alpha+\delta$};
\draw[arrows={-angle 90}, shorten >=9, shorten <=9]  (8,2) -- (8,1);
è\draw[arrows={-angle 90}, shorten >=15, shorten <=15]  (8,2) -- (6,1);
\draw[arrows={-angle 90}, shorten >=15, shorten <=15]  (6,2) -- (8,1);
\draw[arrows={-angle 90}, shorten >=9, shorten <=9]  (6,2) -- (6,1);
\draw[arrows={-angle 90}, shorten >=12, shorten <=12]  (8,1) -- (7,0);
\draw[arrows={-angle 90}, shorten >=12, shorten <=12]  (6,1) -- (7,0);          
        \end{tikzpicture}
  % \end{center}
         %\begin{center}
          
           \vspace{3mm}
    \hspace{3.8cm}
        \begin{tikzpicture}[scale=0.8]
\node at (-2.2,1) {$\xi_{(2,2)}=$};
\node at (0,0) {$1$}; 
\node at (1,1) {$1$};
\node at (-1,1) {$1$};
\node at (1,2) {$1$};
\node at (-1,2) {$1$};

\draw[arrows={-angle 90}, shorten >=9, shorten <=9]  (1,2) -- (1,1);
è\draw[arrows={-angle 90}, shorten >=15, shorten <=15]  (1,2) -- (-1,1);
\draw[arrows={-angle 90}, shorten >=15, shorten <=15]  (-1,2) -- (1,1);
\draw[arrows={-angle 90}, shorten >=9, shorten <=9]  (-1,2) -- (-1,1);
\draw[arrows={-angle 90}, shorten >=12, shorten <=12]  (1,1) -- (0,0);
\draw[arrows={-angle 90}, shorten >=12, shorten <=12]  (-1,1) -- (0,0);           
        \end{tikzpicture}
     %   \end{center}

Above the equivarant cohomology class $(z_{\ell_\bullet})_{\ell_\bullet\in\mathcal{C}(k,n,\omega)}$ is represented by a collection of polynomials arranged on the corresponding vertices of the moment graph (with $z_{\ell_\bullet}$ on the vertex $\ell_\bullet$).
\end{example}

\subsection{Cyclic group action on equivariant cohomology}

Consider the $\mathbb{Z}_n$ action on $T$ given by
\[
m\cdot (\gamma_0, \gamma_1, \ldots, \gamma_n)=(\gamma_0, \gamma_{1+m}, \ldots, \gamma_{n+m}), \quad (m\in\mathbb{Z}_n, \gamma\in T).
\]
Such an action induces a $\mathbb{Z}_n$ action on the character lattice $\mathfrak{X}^*(T)$ via $m\cdot \alpha(\gamma)=\alpha((-m)\cdot \gamma)$ for any $m\in\mathbb{Z}_n, \ \alpha\in \mathfrak{X}^*(T), \ \gamma\in T$. In this way we also get a $\mathbb{Z}_n$ action on $R$, uniquely determined by $m(\epsilon_i)=\epsilon_{i-m}$ and $m(\delta)=\delta$ for any $m,i\in\mathbb{Z}_n$. 

Clearly, the set $\mathcal{C}_{k,n,\omega}$ is also equipped with a $\mathbb{Z}_n$ action given by $m\cdot(\ell_j)_{j\in\mathbb{Z}_n}=(\ell'_j)_{j\in\mathbb{Z}_n}$, where $\ell'_j=\ell_{j-m}$.
\begin{prop}
    There is an action of $\mathbb{Z}_n$ on $H_T^\bullet(X(k,n,\omega),\bQ)$ given by
    \[
    m\cdot \big(z_{\mathcal{\ell_\bullet}}\big)_{\ell_\bullet \in\mathcal{C}({k,n,\omega})}=\big(z'_{\mathcal{\ell_\bullet}}\big)_{\ell_\bullet \in\mathcal{C}({k,n,\omega})}, 
    \]
    where $m\in\mathbb{Z}_n$, $(z_{\mathcal{\ell_\bullet}})\in H_T^\bullet(X(k,n,\omega),\bQ)$, and $z'_{\mathcal{\ell_\bullet}}=m(z_{m\cdot\mathcal{\ell_\bullet}})$.
\end{prop}
\begin{proof}Let $(z_{\mathcal{\ell_\bullet}})\in H_T^\bullet(X(k,n,\omega),\bQ)$ and $m\in\mathbb{Z}_n$. We have to check that for any $\ell_\bullet\in \mathcal{C}({k,n,\omega})$ and any triple $i,j,r$ such that $f_{i,j,r}(\ell_\bullet)$ is well defined, the following holds
\[
mz_{m\cdot\mathcal{\ell_\bullet}}\equiv mz_{m\cdot f_{i,j,r}\mathcal{\ell_\bullet}}\mod \alpha_{\mathcal{\ell_\bullet}, f_{i,j,r}\mathcal{\ell_\bullet}}=\epsilon_j-\epsilon_i+(\ell_i-\ell_j-r)\delta.
\]
Observe that 
\[
m\cdot f_{i,j,r}(\ell_\bullet)_s=
\begin{cases} \ell_{s-m} & s \notin \{i+m,j+m\}\\ \ell_{s-m}-r & s=i+m\\  \ell_{s-m}+r & s=j+m\end{cases}= f_{i+m,j+m,r}(m\cdot \ell_\bullet)_s.
\]
Thus, since $(z_{\mathcal{\ell_\bullet}})\in H_T^\bullet(X(k,n,\omega),\bQ)$, we have
\[
z_{m\cdot\mathcal{\ell_\bullet}}\equiv z_{m\cdot f_{i,j,r}\mathcal{\ell_\bullet}}\mod \epsilon_{j+m}-\epsilon_{i+m}+(\ell_i-\ell_j-r)\delta.
\]
We conclude that 
\[mz_{m\cdot\mathcal{\ell_\bullet}}\equiv mz_{m\cdot f_{i,j,r}\mathcal{\ell_\bullet}}\mod m(\epsilon_{j+m}-\epsilon_{i+m}+(\ell_i-\ell_j-r)\delta)=\epsilon_j-\epsilon_i+(\ell_i-\ell_j-r)\delta.\]
\end{proof}
\begin{example}
    In terms of the KT classes from Example \ref{ex:KT-classes}, the $\mathbb{Z}_2$ action is uniquely determined by 
\[
\sigma\cdot \xi_{(2,2)}=\xi_{(2,2)}, \ \sigma\cdot\xi_{(3,1)}=\xi_{(1,3)}, \ \sigma\cdot \xi_{(4,0)}=\xi_{(0,4)},
\]
where $\sigma$ is the generator of $\mathbb{Z}_2$.

Let $P$ denote $\mathbb{Z}_2$ representation given by the space $R$ equipped with the action introduced above.
As a $\mathbb{Z}_2$ representation, the equivariant cohomology decomposes as
\begin{align*}
H_T^\bullet(X(2,1,2))&\simeq P\xi_{(2,2)}\oplus P(\xi_{(3,1)}+\xi_{(1,3)})\oplus P(\xi_{(4,0)}+\xi_{(0,4)})\\
&\quad\oplus P(\xi_{(3,1)}-\xi_{(1,3)})\oplus P(\xi_{(4,0)}-\xi_{(0,4)}) \\
&\simeq {\bf 1}^0\otimes P\oplus {\bf 1}^2\otimes P\oplus {\bf 1}^4\otimes P\oplus {\bf\epsilon}^2\otimes P\oplus {\bf \epsilon}^4\otimes P,
\end{align*}
where ${\bf 1}^j$ and  ${\bf \epsilon}^j$ represent the one dimensional representation concentrated in degree $j$, on which $1\in\mathbb{Z}_2$ acts via multiplication by 1 and -1, respectively (we recall that the cohomological degree of a cocharacter is 2).
\end{example}
\begin{rem}
It would be interesting to investigate the structure of $\mathbb{Z}_n$ representation on $H_T^\bullet(X(k,n,\omega),\bQ)$ in general.
%, that is how it decomposes into direct sums of irreducible (graded) $\mathbf{Z}_n$-representations.  
%, which we expect to be related to a cyclic sieving phenomenon.    
\end{rem}

\section{Poset structures on the set of fixed points}\label{sec:Poset-structures}
\subsection{Affine permutations, Grassmann necklaces and juggling patterns}\label{sec:affine-permutations} 
In this subsection we briefly recall several combinatorial objects playing an important role in the theory of tnn Grassmannians. The details can be found in \cite{KLS13,Lam16,Rie99,W05}. 

Let ${\mathcal S}_{k,n}$ be the set of $(k,n)$ affine permutations, i.e.  ${\mathcal S}_{k,n}$ consists of bijections $f:\bZ\to\bZ$ such that 
$f(i+n)=f(i)+n$ for all $i\in\bZ$ and $\sum_{i=1}^n (f(i)-i)=kn$. We denote by ${\rm id}_k\in {\mathcal S}_{k,n}$ the permutation given by ${\rm id}_k(j)=j+k$ for all $j\in\bZ$. In particular, $(0,n)$ affine permutations form a group $W_n$ isomorphic to the Weyl group of the affine type $A_{n-1}^{(1)}$.
The group $W_n$ acts on ${\mathcal S}_{k,n}$ by left multiplication. In particular, the map $w\mapsto w{\rm id}_k$ gives a bijection $W_n\to {\mathcal S}_{k,n}$.
Hence the Bruhat order on $W_n$ induces an order on the set of affine permutations.   

A $(k,n)$ affine permutation is called bounded if $i\le f(i)\le i+n$ %f(i+n)$
for all $i\in\bZ$. The set of bounded $(k,n)$ affine permutations is denoted by ${\mathcal B}_{k,n}$. 
The following fact will be important for us:
\begin{equation}
{\mathcal B}_{k,n}  \text{ is a lower order ideal in } {\mathcal S}_{k,n}\simeq W_n. 
\end{equation}

%\subsection{Grassmann necklaces} 
A collection $\mathcal{I}=(I_a)_{a\in[n]}$ of subsets of the set $[n]$ is called a $(k,n)$ Grassmann necklace if $|I_a|=k$ for all $a$ and $I_a\subset I_{a+1}\cup\{a\}$ for all $a\in [n]$ (for $a=n$ we put $I_{a+1}=I_1$). We denote the set of $(k,n)$ Grassmann necklaces by $\mathcal{GN}_{k,n}$. The bijection between $\mathcal{GN}_{k,n}$ and  
${\mathcal B}_{k,n}$ is given by the following rule. For $\mathcal{I}\in \mathcal{GN}_{k,n}$ the corresponding $f\in {\mathcal B}_{k,n}$  fixes all $a$ such that $a\notin I_a$. If $a\in I_a$, then $I_{a+1}=I_a\setminus \{a\}\cup\{b\}$; we put $f(a)=c$, where $a<c\le a+n$ and $b=c$ modulo $n$.  

\begin{example}\label{fI}
For a subset $I\subset [n]$, $|I|=k$ we have a Grassmann necklace ${\mathcal I}=(I,\dots,I)$.  The bounded $(k,n)$ affine permutation $f_I$ corresponding to ${\mathcal I}$ is determined by 
$f_I(a)=\begin{cases}a, & a\notin I,\\ a+n, & a\in I\end{cases}$.  

We denote by $w_I\in W_n$ the Weyl group element corresponding to $f_I$ under the identification $W_n\simeq {\mathcal S}_{k,n}$. 
\end{example}

%\subsection{Juggling patterns}\label{subsec:Jug}
Juggling patterns are close cousins of the Grassmann necklaces. By definition, a collection ${\mathcal J}=(J_1,\dots,J_n)$ of $k$-element subsets of $[n]$  is a $(k,n)$ juggling pattern if 
$\tau_1(J_a\setminus\{n\})\subset J_{a+1}$ for all $a\in [n]$, where $\tau_1(x)=x+1$. Let ${\mathcal Jug}(k,n)$ be the set of $(k,n)$ juggling patterns. 
Then there is a bijection between ${\mathcal Jug}(k,n)$ and  ${\mathcal B}_{k,n}$ given by the following rule. For a juggling pattern ${\mathcal J}=(J_a)_{a=1}^n$ the corresponding element  $f_{\mathcal J}\in {\mathcal B}(k,n)$ is defined by
\begin{equation}\label{JugtoB}
f_{\mathcal J}(a)=\begin{cases} a, & n\notin J_a,\\
n+a+1-x, & n\in J_a, J_{a+1}=\tau_1(J_a\setminus \{n\})\cup \{x\}.
\end{cases}
\end{equation}

\subsection{Generalized bounded affine permutations}\label{sec:bound-affine-permutation}%----------------------------------
As shown in \cite[Section~1]{FLP21}, there is a bijection between cells of $X(k,n,1)$, juggling patterns and bounded affine permutations. Now, we introduce generalized bounded affine permutations which are in bijection with generalized juggling patterns. This gives rise to an alternative way to parameterize the $T$-fixed points of $X(k,n,\omega)$. 
\begin{dfn}\label{def:bounded-affine-permutation}
For $k,n,\omega \in \bN$ with $k \leq n$, a $(k,n,\omega)$ bounded affine permutation is a bijection 
$f:\bZ\to\bZ$ satisfying the following properties:
\begin{enumerate}
    \item $f(i+n)=f(i)+n$ for all $i\in\bZ$,
    \item $\sum_{i=1}^n (f(i)-i)=kn\omega$,
    \item $i \leq f(i) \leq i+\omega n$ for all $i \in \bZ.$
\end{enumerate}    
The set of $(k,n,\omega)$ bounded affine permutations is denoted by $\mathcal{B}_{k,n,\omega}$.
\end{dfn}
Here, Condition (1) is the same as for the $(k,n)$ bounded affine permutations from Section~\ref{sec:affine-permutations}. For $\omega =1$ Condition (2) and (3) are the same as in Section~\ref{sec:affine-permutations}. Without condition (3) we say $f$ is a $(k\cdot \omega,n)$ affine permutation. This definition is valid for all $q \in \bZ$, not only $q =\omega k$. By $\mathcal{S}_{q,n}$ we denote the set of all $(q,n)$ affine permutations.

There is a special $(q,n)$ affine permutation ${\rm id}_{q}$ given by ${\rm id}_k(i)=i+q$. In the setting $q=k\omega$ this is a $(k,n,\omega)$ bounded affine permutation. Without dependence on $\omega$, the length of an affine permutation is defined as 
\[
l(f)=|\{(i,j)\in[n]\times\bZ:\ i<j \text{ and } f(i)>f(j)\}|.
\]
We note that the set of $(0,n)$ affine permutations is a group isomorphic to the affine Weyl group $W_n$ of type $A_{n-1}^{(1)}$.
For general $q$ the group $W_n$ acts freely and transitively on $\mathcal{S}_{q,n}$, because we can write $q=rn+k$ with $0\leq k <n$ and use the same arguments as for $\mathcal{S}_{k,n}$. 
The action of the permutation 
$s_i=(i,i+1)\in W_n$, for $i=0,\dots,n-1$ permutes the values $f(i+rn)$ and $f(i+rn+1)$ for all $r\in\bZ$. This allows to identify
the set of $\mathcal{S}_{q,n}$ with $W_n$ by sending $w\in W_n$ to $w.{\rm id}_{q}$. Hence we obtain an induced order $\le$
on the set $\mathcal{S}_{q,n}$ coming from the Bruhat order on $W_n$. Thus, the unique minimal element is 
${\rm id}_{q}$. 

It is shown in \cite[Lemma~3.6]{KLS13} that ${\mathcal B}_{k,n,1}$
is a lower order ideal in ${\mathcal S}_{k,n}$. With the same arguments it follows that ${\mathcal B}_{k,n,\omega}$
is a lower order ideal in ${\mathcal S}_{k\cdot \omega,n} \cong W_n$. For $f,g\in {\mathcal B}_{k,n,\omega}$ we write $f\le_{\mathcal B} g$ for the order induced by the Bruhat order on $W_n$.

By \cite[Theorem 6.2]{Lam16}, there is an order preserving bijection between the set ${\mathcal B}_{k,n,1}$ and the set ${\mathcal Jug}(k,n,1)$. Before generalizing this to arbitrary $\omega$, we introduce an alternative parametrization of the $T$-fixed points of $X(k,n,\omega)$ which is closer to the definition of bounded affine permutations.
\begin{prop}\label{prop:cell-parametrization}
For $k,n,\omega \in \bN$ with $k \leq n$, there is a bijection between ${\mathcal Jug}(k,n,\omega)$ and 
\[\Big\{ (\ell_j)_{j \in \bZ_n} \in [0,\omega n]^{\bZ_n}  \ :  \ \sum_{j \in \bZ_n}\ell_j = kn\omega, \ \ \big\{ j+ \ell_j\,\mathrm{mod}\, n : j \in \bZ_n\big\} = \bZ_n\Big\}. \]
\end{prop}
\begin{proof}
By Proposition~\ref{prop:bijection_jug-cells} we have to show that for $\ell_\bullet\in [0,\omega n]^{\bZ_n}$ the following are equivalent:
\begin{enumerate}
    \item $\bdim \, U_{\ell_\bullet} = (k\omega, \dots, k\omega)\in \bN^{\bZ_n}$,
    \item \[\sum_{j \in \bZ_n}\ell_j = kn\omega \quad \mr{and} \quad \big\{ j+ \ell_j \mod n : j \in \bZ_n\big\} = \bZ_n.\]
\end{enumerate}
The length tuple $z_\bullet$ with $z_j=k\omega$ for all $j \in \bZ_n$ satisfies (1) and (2). It represents the unique zero-dimensional cell $C_0$ of $X(k,n,\omega)$ and $C_0$ is contained in the closure of every other cell of $X(k,n,\omega)$. Hence, by \cite[Theorem~6.15]{LaPu2020} every other $\ell_\bullet \in \mathcal{C}_{k,n,\omega}$ is obtained from $z_\bullet$ by a sequence of cut and paste moves $f_{i,j,r} : \mathcal{C}_{k,n,\omega} \to \mathcal{C}_{k,n,\omega}$. These moves preserve the properties (1) and (2). Hence (1) implies (2).

The second part of (2) implies that  $\bdim \, U_{\ell_\bullet} = (L/n, \dots, L/n)$ where $L = \sum_{j \in \bZ_n}\ell_j$. Together with the first part of (2) this implies (1).
\end{proof}
\begin{lem}\label{lem:cell-parametrisations}
For $k,n,\omega \in \bN$ with $k \leq n$, there are bijections between ${\mathcal Jug}(k,n,\omega)$, ${\mathcal C}(k,n,\omega)$ and ${\mathcal B}_{k,n,\omega}$.   
\end{lem}
\begin{proof}
Proposition~\ref{prop:bijection_jug-cells} gives the first bijection. We define the map $\psi: {\mathcal C}(k,n,\omega) \to {\mathcal B}_{k,n,\omega}$ sending $\ell_\bullet$ to the map $f: \bZ \to \bZ$ with 
$f(j'):=j'+\ell_j$ for all $j' \in \bZ$ with $j'=j \mod n$.
It follows from the second parametrization of cells from Proposition~\ref{prop:cell-parametrization} that this $f$ is bijective and satisfies part (2) of Definition~\ref{def:bounded-affine-permutation}. $\ell_j \in [0,\omega n]$ implies part (3) of that definition. The inverse map $\psi^{-1}$ sends $f$ to $(f(j)-j)_{j\in\bZ_n}$.
\end{proof}
\subsection{Partial orders on the set of cells}\label{sec:partial-orders}%-----------------------------
In this section we introduce partial orders on the sets ${\mathcal Jug}(k,n,\omega)$, ${\mathcal C}(k,n,\omega)$ and ${\mathcal B}_{k,n,\omega}$ and examine how they are related under the bijections from Lemma~\ref{lem:cell-parametrisations}.

For $\mathcal{J}_\bullet, \mathcal{J}_\bullet' \in {\mathcal Jug}(k,n,\omega)$ we write $\mathcal{J}_\bullet \geq_{\mathcal J} \mathcal{J}_\bullet'$ iff $j^{(i)}_r \leq j'^{(i)}_r$ for all $i \in \bZ_n$ and $r \in[k\omega]$ where we order each $J_i \in \binom{[n\omega]}{k\omega}$ as
\[ \big(j^{(i)}_1 <  j^{(i)}_2 < \ldots < j^{(i)}_{k\omega} \big). \]
Given two elements $\ell_\bullet, \ell_\bullet' \in {\mathcal C}(k,n,\omega)$ we write $\ell_\bullet \geq_{\mathcal C} \ell_\bullet'$ if there exists an oriented path from $\ell_\bullet$ to $\ell_\bullet'$ in the moment graph for the $T$-action on $X(k,n,\omega)$. For $p,p' \in X(k,n,\omega)^T$ we write $p' \preceq p$ if $\overline{ C_p}$ contains $p'$. By Theorem~\ref{trm:cells-are-strata} we obtain the same partial order $\preceq$ if we consider closures of $\mr{Aut}_{\Delta_n}(U_{\omega n})$-orbits or strata. Recall that the partial order $\geq_{\mathcal B}$ on ${\mathcal B}_{k,n,\omega}$ is induced by the Bruhat order on $W_n$.
\begin{example}
For the $T$-fixed points from Example~\ref{ex:T-fixed-points}, the bijections from Lemma~\ref{lem:cell-parametrisations} are as follows:
\begin{align*}
    (2,2) & \longleftrightarrow \big( \{3,4 \}, \{3,4 \}\big) \longleftrightarrow f \ \ \mathrm{with} \ f(1):= 1+2, f(2):=2+2,\\
    (3,1) & \longleftrightarrow \big( \{2,4 \}, \{3,4 \}\big) \longleftrightarrow f \ \ \mathrm{with} \ f(1):= 1+3, f(2):=2+1,\\
    (1,3) & \longleftrightarrow \big( \{3,4 \}, \{2,4 \}\big) \longleftrightarrow f \ \ \mathrm{with} \ f(1):= 1+1, f(2):=2+3,\\
    (4,0) & \longleftrightarrow \big( \{2,4 \}, \{1,3 \}\big) \longleftrightarrow f \ \ \mathrm{with} \ f(1):= 1+4, f(2):=2+0,\\
    (0,4) & \longleftrightarrow \big( \{1,3 \}, \{2,4 \}\big) \longleftrightarrow f \ \ \mathrm{with} \ f(1):= 1+0, f(2):=2+4.
\end{align*}
and all poset structures are the same as the one induced by Example~\ref{ex:one-dim-orbits}.
\end{example}
The next theorem shows that this identification of poset structures was no coincidence.
\begin{thm}\label{thm:poset-iso}
For $k,n,\omega \in \bN$ with $k \leq n$, there are order preserving poset isomorphisms  between ${\mathcal Jug}(k,n,\omega)$, ${\mathcal C}(k,n,\omega)$,  $X(k,n,\omega)^T$ and ${\mathcal B}_{k,n,\omega}$. 
\end{thm}
\begin{cor}\label{cor:about-cells-and-their-closures}
The closure of every cell in $X(k,n,\omega)$ is obtained as
\[ \overline{ C_\mathcal{J}} = \bigcup_{\mathcal{J}' \in {\mathcal Jug}(k,n,\omega) \ \mathrm{s.t.:} \ \mathcal{J}' \leq_{\mathcal \mathcal{J}} \mathcal{J}} C_{\mathcal{J}'}.\]
Moreover the moment graph of $\overline{ C_\mathcal{J}}$ is the full subgraph of the graph described in Lemma~\ref{lem:moment-graph-edges} on the vertices corresponding to $\mathcal{J}' \leq_{\mathcal J} \mathcal{J}$. The dimension of $C_\mathcal{J}$ is the number of edges in the moment graph starting at $\mathcal{J}$. This equals the length of the corresponding bounded affine permutation.
\end{cor}
\begin{proof}
The description of the closure is obtained from Theorem~\ref{thm:poset-iso} in combination with Remark~\ref{rem:cell-closures}. Hence the moment graph of the cell closure is the full subgraph on the vertices which are smaller with respect to any of the partial orders. 
Finally the dimension formula is obtained from the embedding into the affine flag variety as described in Section~\ref{sec:quiver-Grass-in-affine-flag}. 
\end{proof}
Below, we relate the poset structures of the cells and generalized juggling patterns. The relation to affine permutations is examined in Section~\ref{sec:quiver-Grass-in-affine-flag}.
\begin{proof}[Proof of Theorem~\ref{thm:poset-iso} ($ \leq_{\mathcal J} \iff \leq_{\mathcal C} \iff \preceq$)]
The isomorphisms on the level of sets were introduced in Lemma~\ref{lem:cell-parametrisations} and Lemma~\ref{lem:T-fixed-point-parametrisation}. It remains to show that they preserve the poset structures. The partial order  $\leq_{\mathcal C}$ on ${\mathcal C}(k,n,\omega)$ is obtained from the edges in the moment graph as described in Lemma~\ref{lem:moment-graph-edges}. The property $\ell_i > \ell_j+r$  of the cut and paste moves implies that $\varphi^{-1}(\ell_\bullet)\geq_{\mathcal{J}} \varphi^{-1}(f_{i,j,k}(\ell_\bullet))$ where $\varphi$ is the map from Proposition~\ref{prop:bijection_jug-cells}. Hence $\ell_\bullet \geq_{\mathcal{C}} \ell_\bullet'$ implies $\varphi^{-1}(\ell_\bullet) \geq_{\mathcal{J}} \varphi^{-1}(\ell_\bullet')$ for all $\ell_\bullet,\ell_\bullet' \in {\mathcal C}(k,n,\omega)$. 

Starting with $\mathcal{J}_\bullet, \mathcal{J}_\bullet' \in {\mathcal Jug}(k,n,\omega)$ such that $\mathcal{J}_\bullet  \geq_{\mathcal{\mathcal{J}}} \mathcal{J}_\bullet'$ we construct a path from $\ell_\bullet:=\varphi(\mathcal{J}_\bullet)$ to $\ell_\bullet':=\varphi(\mathcal{J}_\bullet')$ in the moment graph inductively:  Let $d:=\# \{ i \in \bZ_n : J_i \neq J_i'\}$. For $d=0$ both juggling patterns are equal and there is nothing to show. If $d>0$ there exists an $s \in \bZ_n$ such that $ J_s \neq J_s'$ and an $r \in [\omega k]$ such that $j^{(s)}_p=j'^{(s)}_p$ for all $p \in [r+1,\omega k]$ and $j^{(s)}_r<j'^{(s)}_r$. Here we assume that the sets belonging to the juggling patterns are ordered increasingly as introduced in the beginning of Section~\ref{sec:partial-orders}. Now, we determine a path $\ell_\bullet \to \ell_\bullet''$ in the moment graph with $\varphi^{-1}(\ell_\bullet'') =: \mathcal{J}_\bullet'' \geq_{\mathcal{J}} \mathcal{J}_\bullet'$ and $d>\# \{ i \in \bZ_n : J_i'' \neq J_i'\}$. 

The points $j^{(s)}_r$ and $j'^{(s)}_r$ live on two different indecomposable summands of the representation $U_{\omega n}$ indexed by $a,b \in \bZ_n$ which are obtained from $\mathcal{J}_\bullet, \mathcal{J}_\bullet'$ by the map $\varphi$ as described in the proof of Proposition~\ref{prop:bijection_jug-cells}. This gives rise to the cut and paste map $f_{a,b,r}$ where $r:=b-a+\ell_a-\ell_b$. By construction it follows that $\ell_\bullet >_{\mathcal{C}} f_{a,b,r}(\ell_\bullet)$ and $\varphi^{-1}(f_{a,b,r}(\ell_\bullet)) \geq_{\mathcal{J}} \mathcal{J}_\bullet'$. We can apply this construction of an edge in the moment graph recursively until we reach a point $\ell_\bullet''$ with $\varphi^{-1}(\ell_\bullet'') =: \mathcal{J}_\bullet'' \geq_{\mathcal{J}} \mathcal{J}_\bullet'$ and $J_s''=J_s'$. This implies $d>\# \{ i \in \bZ_n : J_i'' \neq J_i'\}$ and finishes the inductive step. Hence the partial orders $ \leq_{\mathcal J} $ and $\leq_{\mathcal C}$ are equivalent.

With the explicit description of the cells as attracting sets of the fixed points, it is straightforward to check that for $p,p' \in X(k,n,\omega)^T$ the moment graph contains a path from $p$ to $p'$ if and only if $\overline{ C_p}$ contains $p'$. This implies the equivalence of the partial orders $\leq_{\mathcal C}$ and $\preceq$.
\end{proof}
\subsection{Poincar\'e polynomials}%-----------------------------
The Poincar\'e polynomial of the quiver Grassmannian $X(k,n,\omega)$ is obtained as 
\[ P_{k,n,\omega}(q) = \sum_{p \in X(k,n,\omega)^T} q^{\dim_\bC C_p}.\]
From the computations in the previous section we obtain the following formula
\begin{lem}
For $k,n,\omega \in \bN$ with $k \leq n$, the Poincar\'e polynomial of $X(k,n,\omega)$ is
\[ P_{k,n,\omega}(q) = \sum_{f \in \mathcal{B}_{k,n,\omega}} q^{l(f)}.\]
Here $l(f)$ denotes the length of the bounded affine permutation $f$ as defined in Section~\ref{sec:bound-affine-permutation}.
\end{lem}
\begin{proof}
This follows immediately from Corollary~\ref{cor:about-cells-and-their-closures}.
\end{proof}

%----------------------------------------------------------------------
\section{Affine flag varieties}\label{sec:Affine-flag}%----------------------------------
%----------------------------------------------------------------------

In this section we recall some basics from the theory of affine flag varieties in type $A$ and Sato Grassmannians (\cite{Kum02,FFR17,Pue2020}). We use this material in the next section for the explicitl construction of the embeddings of the quiver Grassmannians $X(k,n,\omega)$ into the affine flag varieties.

\subsection{Notation}
Let $\widehat{\msl_n}$ denote the affine Kac-Moody Lie algebra of type $A_{n-1}^{(1)}$. Explicitly, 
$\widehat{\msl_n}=\msl_n\otimes \bC[t,t^{-1}]\oplus\bC K\oplus\bC d$, where $K$ is  central and $d$ is the derivation.
Let us fix a Cartan decomposition $\msl_n=\fn\oplus \fh\oplus\fn_-$ and the Borel
subalgebra $\fb=\fn\oplus\fh$. Then the Iwahori subalgebra of $\widehat{\msl_n}$ is
given by $\fb\otimes 1\oplus\msl_n\otimes t\bC[t]$.

We denote the Weyl group of $\widehat{\msl_n}$ by $W_n$. 
For $n=2$ the Weyl group is generated by $s_0$ and $s_1$ subject to the relations
$s_0^2=s_1^2=e$. For $n>2$ the group $W_n$ is generated by
reflections $s_0,s_1,\dots,s_{n-1}$ subject to the defining relations
\begin{align*}
s_i^2=e,\ i=0,\dots,n-1,\quad s_is_j=s_js_i,\ |i-j|>1,\\
s_is_{i+1}s_i=s_{i+1}s_is_{i+1},\ i=0,\dots, n-1
\end{align*}
(here and below we set $s_n=s_0$). The group $W_n$ can be also realized as the group of $(0,n)$ affine permutations, i.e. bijections $f:\bZ\to\bZ$ subject to the conditions $f(i+n)=f(i)+n$ for all $i$ and $\sum_{i=1}^n (f(i)-i)=0$.

Let $\widehat{SL_n}$ be the affine group with the Lie algebra $\widehat{\msl_n}$.
This group contains the finite torus $\exp(\fh)$ and the two-dimensional torus 
$(\bC^*)^2=\exp(\bC K \oplus\bC d)$.   
We denote by $P_i\subset \widehat{SL_n}$, $i=0,\dots,n-1$  the maximal parabolic subgroups. Then the affine Grassmannians are defined as the quotients
$\widehat{SL_n}/P_i$. 

Let ${\mathfrak B}\subset \widehat{SL_n}$  be the Iwahori subgroup. More precisely, ${\mathfrak B}$
consists of matrices $A(t)\in SL_n(\bC[t])\subset \widehat{SL_n}$ such that $A(0)$ is upper triangular. In particular, the Lie algebra of ${\mathfrak B}$ is $\fb \otimes 1 \oplus \msl_n\otimes t\bC[t]$ 

Let $\eF_n\simeq \widehat{SL_n}/{\mathfrak B}$ be the affine flag variety for the group $\widehat{SL_n}$. One has the natural embedding of the affine flag variety into the product of affine Grassmannians $\eF_n\subset \prod_{i=0}^{n-1} \widehat{SL_n}/P_i$. 
We note that  the fixed points of $\eF_n$ with respect to the torus of $\widehat{SL_n}$  are labeled by the elements of the group $W_n$. For $w\in W_n$
let $p_w\in \eF_n$ be the corresponding torus fixed point. 

The affine flag variety is an ind-variety, i.e. the inductive limit of finite-dimensional projective algebraic varieties. Namely, for an element 
$w\in W_n$ we denote the corresponding (finite-dimensional) affine Schubert variety $\overline{{\mathfrak B}.p_w}$ by $X_w$. 
Then  $\eF_n=\bigcup_{w\in W_n} X_w$.

\subsection{Sato Grassmannians}
The affine Grassmannians enjoy  explicit embeddings to the Sato Grassmannians. Hence the affine flag variety can be realized inside the product of Sato Grassmannians. We provide some details below. 

The Sato Grassmannian ${\rm SGr}^{(i)}$, $i\in\bZ$ consists of subspaces $V\subset \bC[t,t^{-1}]$ such that
\begin{itemize}
    \item $t^N\bC[t^{-1}] \supset V\supset t^{-N}\bC[t^{-1}]$ for some $N\in\bZ$,
    \item $\dim V/t^{-N}\bC[t^{-1}]=i+N$.
\end{itemize}
For example, the subspace $\accentset{\circ}{V}^{(i)}={\rm span}\{t^j:\ j\le i\}$ belongs to ${\rm SGr}^{(i)}$.

Sato Grassmannians are ind varieties (in particular, they can be realized as inductive limits of finite-dimensional Grassmann varieties). They enjoy a Pl\"ucker embedding
into the projective space $F=\Lambda^{\infty/2}(\bC[t,t^{-1}])$ of semi-infinite forms. 
The space $F$ is spanned by infinite wedge products   
\[
t^L = t^{l_1} \wedge t^{l_2} \wedge \dots , \qquad  L= (l_1,l_2,\dots ),   
\]
where $l_1>l_2>\dots$ and $l_{s+1}=l_s-1$ for $s$ large enough. One has the charge decomposition
$F=\bigoplus_{i\in\bZ} F^{(i)}$, where $F^{(i)}$ is spanned by wedges $t^L$ such that 
$l_s=i-s+1$ for $s$ large enough. Then one has the Pl\"ucker embedding 
${\rm SGr}^{(i)}\hookrightarrow \bP(F^{(i)})$.

\begin{rem}
Let $|i\rangle\in F^{(i)}$ be the charge $i$ vacuum vector, explicitly given by 
$|i\rangle = t^i \wedge t^{i-1}\wedge \dots$. Then the image of the space $\accentset{\circ}{V}^{(i)}\in {\rm SGr}^{(i)}$ inside
$\bP(F^{(i)})$ coincides with the line containing $|i\rangle$. We also note that each space
$F^{(i)}$ is endowed with the action of the infinite-dimensional Heisenberg algebra.
As a module over the Heisenberg algebra $F^{(i)}$ is isomorphic to a Fock module.
\end{rem}

The affine flag variety $\eF_n$ is realized inside the product $\prod_{i\in \bZ} {\rm SGr}^{(i)}$ as the set of collections $(V_i)_{i\in\bZ}$ such that 
\begin{itemize}
    \item $V_i\subset V_{i+1}$,
    \item $V_{i+n}=t^n V_i$.
\end{itemize}
In particular, the collection $(\accentset{\circ}{V}^{(i)})_i$ corresponds to the coset of the identity
in $\eF_n\simeq \widehat{SL_n}/{\mathfrak B}$.

\begin{rem}
\label{wexplicit}
Given an element $w\in W_n$ the corresponding torus fixed point $p_w\in\eF_n$ is given by 
\[
\prod_{i\in\bZ} {\mathrm{ span}}(t^{w(j)}, j\le i) \in\prod_{i\in\bZ} {\rm SGr}^{(i)}.
\]
\end{rem}

%----------------------------------------------------------------------
\section{Quiver Grassmannians inside affine flags}\label{sec:quiver-Grass-in-affine-flag}

\bigskip

%----------------------------------------------------------------------
The goal of this section is to construct explicitly the embedding of the quiver Grassmannians $X(k,n,\omega)$ into the affine flag varieties. The 
image is described as a union of certain Schubert varieties. We note that
due to the identification of our quiver Grassmannians with the local models
of Shimura varieties %(see Appendix \ref{local}) 
(c.f. Remark~\ref{rem:Shimura-var}) such an embedding can be
found e.g. in \cite{HG02,PRS13}. Our goal in this section is to construct the 
embedding explicitly in the language of quiver Grassmannians and Sato Grassmannians
and to label the Schubert varieties showing up in the image via generalized 
juggling patterns and bounded affine permutations.

\subsection{The construction}
By definition, the quiver Grassmannian $X(k,n,\omega)$ sits inside the product  
$\prod_{i=0}^{n-1} {\rm Gr}_{k\omega}(M^{(i)})$, where $M^{(i)}$ are $n\omega$-dimensional vector spaces with bases $v^{(i)}_j$, $j\in [n\omega]$ (here we
identify $\bZ/n\bZ$ with the set $0,\dots,n-1$). Using the Pl\"ucker embeddings for ${\rm Gr}_{k\omega}(M^{(i)})$ we obtain the embedding
\begin{equation}\label{omegaP}
X(k,n,\omega)\subset  \prod_{i=0}^{n-1} \bP(\Lambda^{k\omega}(M^{(i)}))
\end{equation}

Now let us construct the embeddings 
$\Psi^{(i)}: {\rm Gr}_{k\omega}(M^{(i)})\hookrightarrow {\rm SGr}^{(i)}$. 
We define the map $\psi^{(i)}: M^{(i)}\to \bC[t,t^{-1}]$ by the formula
\begin{equation}\label{vt}
\psi^{(i)} v^{(i)}_{n\omega+1-j} = t^{i-k\omega+j}, \ j\in [n\omega].
\end{equation}
In particular, the image of $\psi^{(i)}$ is spanned by $t^{i-k\omega+1},\dots, t^{i-k\omega+n\omega}$.
\begin{rem}
Let us briefly explain the shifts of indices showing up in \eqref{vt} (the details 
are given in the proofs below). Recall the map
$\tau_1$ defined by $v^{(i)}_j \mapsto  v^{(i+1)}_{j+1}$. The subscript $n\omega +1-j$ in the
left hand side and the power $i-k\omega+j$ are chosen in such a way that the map
$\tau_1$ induces the identity map when translated to $\bC[t,t^{-1}]$. We also note that the term $-k\omega$ in the power of $t$ in the right hand side guaranties that
the induced map $\Lambda^{k\omega}(M^{(i)})\to F$, defined by 
$$\mathrm{ span}\{w_1,\dots,w_{k\omega}\}\mapsto w_1\wedge \dots \wedge w_{k\omega}\wedge |i-k\omega\rangle,$$ lands in the subspace $F^{(i)}$ of the charge $i$
forms. 
\end{rem}

Now the maps $\Psi^{(i)}$ are defined as follows: 
\begin{equation}\label{Psi}
\Psi^{(i)} (U^{(i)}) = \psi^{(i)} U^{(i)} \oplus {\rm span}\{t^j:\ j\le i-k\omega\}. 
\end{equation}
The following lemma is obvious.
\begin{lem}
The image $\Psi^{(i)} {\rm Gr}_{k\omega}(M^{(i)})$ belongs to ${\rm SGr}^{(i)}$. 
\end{lem}

 Using the maps  $\Psi^{(i)}$ one obtains the embedding  
\begin{equation}\label{A}
\Psi: X(k,n,\omega)\to \prod_{i\in\bZ} {\rm SGr}^{(i)}.
\end{equation}
To be precise, $\Psi(X(k,n,\omega))$ sits
inside the product $\prod_{i=0}^{n-1} {\rm SGr}^{(i)}$. We consider a larger product 
$\prod_{i\in\bZ} {\rm SGr}^{(i)}$
for all integers $i$ via the obvious extension rule $M^{(i+n)}=M^{(i)}$, 
${\rm Gr}_{k\omega}(M^{(i)})={\rm Gr}_{k\omega}(M^{(i+n)})$.
This gives the desired realization \eqref{A}.

\begin{rem}
Let us consider the point $\accentset{\circ} U = (\accentset{\circ} U^{(i)})_{i\in\bZ}$ defined by
$\accentset{\circ} U^{(i)} = {\rm span} \{v^{(i)}_j, \ (n-k)\omega+1\le j\le n\omega\}$. Then clearly 
$\accentset{\circ} U\in X(k,n,\omega)$ and $\Psi\accentset{\circ} U$ is the line passing through 
the product of highest weight vectors $|i\rangle$. 
\end{rem}

Recall the embeddings $\eF_n\subset \prod_{i\in\bZ} {\rm SGr}^{(i)}\subset \prod_{i\in\bZ} \bP(F^{(i)})$ from Section~\ref{sec:Affine-flag}.
In what follows we identify the affine flag variety with its image inside the product of Sato 
Grassmannians.  
\begin{lem}
The image of $\Psi$ belongs to $\eF_n$.
\end{lem}
\begin{proof}
Let $U=(U^{(i)})_{i\in\bZ}$ be a point in $X(k,n,\omega)$ (recall that we use
the extension $U^{(i)}=U^{(i+n)}$ inside $M^{(i)}=M^{(i+n)}$). 
We have to show that
\begin{itemize}
    \item $\Psi^{(i)}U^{(i)} \subset \Psi^{(i+1)}U^{(i+1)}$,
    \item $\Psi^{(i+n)}U^{(i+n)}=t^n \Psi^{(i)}U^{(i)}$.
\end{itemize}
By definition,  
\[
\psi^{(i+n)} = t^n \psi^{(i)} \text{ and } {\rm span}\{t^j:\ j\le i+n-k\omega\} = t^n {\rm span}\{t^j:\ j\le i-k\omega\}. 
\]
Hence the second property $\Psi^{(i+n)}U^{(i+n)}=t^n \Psi^{(i)}U^{(i)}$ holds.

To prove the first property we need to show that
\[
\psi^{(i)} U^{(i)} \oplus {\rm span}\{t^j:\ j\le i-k\omega\} \subset 
\psi^{(i+1)} U^{(i+1)} \oplus {\rm span}\{t^j:\ j\le i+1-k\omega\}.
\]
Clearly it suffices to show that 
\begin{equation}\label{in}
\psi^{(i)} U^{(i)} \subset 
\psi^{(i+1)} U^{(i+1)} \oplus {\rm span}\{t^j:\ j\le i+1-k\omega\}.
\end{equation}
Recall that $\tau_1 U^{(i)}\subset U^{(i+1)}$, where 
$\tau_1 v^{(i)}_j = v^{(i+1)}_{j+1}$ unless $j=n\omega$ and $\tau_1 v^{(i)}_{n\omega}=0$.
For a vector $u=\sum_{j=1}^{n\omega}  r_j v^{(i)}_{n\omega+1-j}\in U^{(i)}$ one has:
\[
\psi^{(i)} u = \sum_{j=1}^{n\omega}  r_j t^{i-k\omega+j},\ \quad
\psi^{(i+1)} (\tau_1u) = \sum_{j=2}^{n\omega}  r_j t^{i-k\omega+j}.
\]
Since $\tau_1u\in U^{(i+1)}$ we obtain \eqref{in}. 
\end{proof}

\subsection{Description of the image}
Now our goal is to identify the image $\Psi X(k,n,\omega)$ inside the affine flag variety. Recall (see Theorems \ref{thm:geometric-properties} 
and \ref{trm:cells-are-strata}) that $X(k,n,\omega)$ has $\binom{n}{k}$ irreducible
components $X_I(k,n,\omega)$ labeled by the cardinality $k$ subsets $I\subset [n]$. Each irreducible component is a closure of a cell $C_I$ containing 
a unique torus fixed point $p_I$. The cell $C_I$ is equal to the  
${\rm Aut}_{\Delta_n}(U_{n\omega})$ orbit passing through $p_I$. Our goal here is two-fold.
First, we show that after the embedding of $X(k,n,\omega)$ into the affine flag variety
the action of the automorphism group  translates into the action of the Iwahori subgroup 
${\mathfrak B}$. Second, we compute the images $\Psi p_I$ of the torus fixed points, i.e.
we find the affine Weyl group elements $w(I)$ such that $\Psi p_I=p_{w(I)}$. As a corollary
we conclude that the embedding $\Psi$ realizes the quiver Grassmannian $X(k,n,\omega)$
as the union of the Schubert varieties $X_{w(I)}$ inside the affine flag variety $\eF_n$.

\begin{lem}
The point $p_I$ is determined by the condition $v^{(i)}_1\in p_I^{(i)}$ for all $i\in I$.
\end{lem}
\begin{proof}
We note that, since $\tau_1 U^{(i)}\subset U^{(i+1)}$ for any point $U\in X(k,n,\omega)$ with $U= (U^{(i)})_{i\in\bZ}$, the condition $v^{(i)}_1\in p_I^{(i)}$ for $i\in I$ says that 
\begin{gather*}
p_I^{(i_0)} \ni 
v^{(i_0)}_{1+(i_0-i)+nr} \text{ for all } i\in I, i\le i_0 \text{ and } r=0,\dots,\omega-1,\\
p_I^{(i_0)} \ni 
v^{(i_0)}_{1+(i-i_0)+n(r+1)} \text{ for all } i\in I, i>i_0 \text{ and } r=0,\dots,\omega-1.
\end{gather*}
Hence for each pair $i\in I$, $r=0,\dots,\omega-1$ we obtain a basis vector in $U^{(i_0)}$.
Since $|I|=k$ and $\dim U^{(i_0)}=k\omega$ we obtain the desired claim. The only thing to
mention is that one easily sees that a point $p_I$ as above does not belong to the closure of the ${\rm Aut}_{\Delta_n}(U_{n\omega})$ orbit of any other torus fixed point (thanks to the
explicit form of the automorphism group elements from Proposition \ref{Aut}). Hence the 
closure of the the orbit ${\rm Aut}_{\Delta_n}(U_{n\omega})p_I$ is an irreducible component for any $I$.
\end{proof}

\begin{rem}
The ($\omega$-generalized) juggling pattern $(J_1,\dots,J_n)$ corresponding to the point $p_I$ is given by 
\begin{align*}
    J_{i_0}=&\{1+(i_0-i)+n(r-1), \ i\in I, i\le i_0, r\in [\omega]\} \\
    &\quad\ \cup 
\{1+(i-i_0)+nr, \ i\in I, i> i_0, r\in [\omega]\}
\end{align*}
for all $i_0\in [n]$.
\end{rem}

For $I\in\binom{[n]}{k}$ we denote by $w(I)\in W_n$ the Weyl group element such that
$\Psi p_I=p_{w(I)}$.
\begin{cor}
The element $w(I)\in W_n$ is defined by 
\[
w(I): i\mapsto  \begin{cases} i-k\omega, & i\notin I,\\ i-k\omega+n\omega, & i\in I \end{cases}
\]
for all $i=1,\dots,n$.
\end{cor}
\begin{proof}
Let $p_I=(U^{(i)})_{i\in [n]}$. Then one has $U^{(i)}={\rm span}\{v^{(i)}_j:\ j\in R_i\}$ for certain
subsets $R_i\subset [\omega n]$. We note that the following holds true:
\begin{equation}\label{RR}
R_i=\begin{cases} \{a+1,\ a\in R_{i-1}\}, & i\notin I,\\ 
\{a+1,\ a\in R_{i-1}, a\ne \omega n\}\cup\{1\}, & i\in I.\end{cases}
\end{equation}

Now let $\Psi p_I =(\Psi p_I)_i$ ($i\in\bZ$). Then each space $(\Psi p_I)_i$ is a vector space
spanned by the elements $t^a$, $a\in S_i$ for a subset $S_i\subset \bZ$. By definition 
$w(I)(i)=S_i\setminus S_{i-1}$ (see Remarks \ref{wexplicit}). 
Formulas \eqref{vt} and \eqref{Psi} say that
\[
S_i = \{i-k\omega+n\omega+1-r, r\in R_i\} \cup \bZ_{\le i-k\omega}.
\]
Taking into account equality \eqref{RR} we obtain
\[
S_i\setminus S_{i-1} = 
\begin{cases} i-k\omega, & i\notin I,\\ i-k\omega+n\omega, & i\in I \end{cases}.
\]
\end{proof}
\begin{rem}
Since $|I|=k$, one has $\sum_{i\in[n]} (w(I)(i)-i)=0$.
\end{rem}

Let ${\widetilde{\mathfrak B}}\subset GL_n(\bC[z])$ be the Iwahori subgroup consisting of such
matrices $\tilde g(z)$ that $\tilde g(0)$ is lower triangular.
\begin{prop}\label{AutB}
The map $\Psi$ translates the action of ${\rm Aut}_{\Delta_n}(U_{n\omega})$ on $X(k,n,\omega)$ to the ${\widetilde{\mathfrak B}}$ action on its image. 
\end{prop}
\begin{proof}
Let us recall the explicit action of the Iwahori group on the affine flag
variety sitting inside the product of Sato Grassmannians. 
Let us fix an $n$-dimensional space $W$ with a basis $w_1,\dots,w_n$.
We identify
the space $\bC[t,t^{-1}]$ with $W\otimes \bC[z,z^{-1}]$ via the map
\[
\varphi:  t^{(n-k)\omega - rn - s} \mapsto w_s\otimes z^r,\ r\in\bZ, s=0,\dots,n-1.
\]
Let us comment on the choice of the map $\varphi$. 
Recall (see \eqref{vt}) the mappings $\psi^{(i)}: M^{(i)}\to \bC[t,t^{-1}]$, sending 
$v^{(i)}_{n\omega+1-j}$ to $t^{i-k\omega+j}$.  
The map $\varphi$ is chosen in such a way that the composition map $\varphi \psi^{(0)} : M^{(0)}\to W\otimes \bC[z,z^{-1}]$ (see \eqref{vt}) is given by
\[
\varphi \psi^{(0)}: v^{(0)}_{rn+s}\mapsto w_s\otimes z^r.
\]
The general composition map $\varphi \psi^{(i)} : M^{(i)}\to W\otimes \bC[z,z^{-1}]$ reads as
\[
\varphi \psi^{(i)}: v^{(i)}_{rn+s+i}\mapsto w_s\otimes z^r.
\]

The natural action of $GL_n(\bC[z])$ (and thus of the  Iwahori subgroup $\widetilde {\mathfrak B}$) on the
space $W\otimes \bC[z,z^{-1}]$ induces the action on $\bC[t,t^{-1}]$ via the identification above.  Each Iwahori group element $\tilde g(z)$ induces the linear endomorphism of  $F^{(i)}$
and a map on $\eF_n$ embedded into the product of the projective spaces  $\bP(F^{(i)}$.

Recall (see Proposition \ref{Aut}) that the elements of the automorphism group ${\rm Aut}_{\Delta_n}(U_{n\omega})$
are the collections of maps $A_i\in {\rm End}(M^{(i)})$ with
\[ 
A_i = \begin{pmatrix}
a^{(i)}_{1,1} & & & &\\
a^{(i)}_{2,1}& a^{(i-1)}_{1,1} & & &\\
\vdots & \vdots & \ddots& &\\
a^{(i)}_{n\omega-1,1} & a^{(i-1)}_{n\omega-2,1} & \hdots & a^{(i-n\omega+2)}_{1,1} &\\
a^{(i)}_{n\omega,1} & a^{(i-1)}_{n\omega-1,1} & \hdots & a^{(i-n\omega-2)}_{2,1} & a^{(i-n\omega+1)}_{1,1}
\end{pmatrix} 
\]
written in the basis $v^{(i)}_j$. 

Given a collection of maps $A=(A_i)\in {\rm Aut}_{\Delta_n}(U_{n\omega})$
we construct an Iwahori group element $\tilde g(z)\in GL_n(\bC[z])$ such that
$\Psi A = \tilde g(z)$.
Recall that the map $\Psi: X(k,n,\omega)\to \eF_n$ is induced by the maps \eqref{vt}.
Then one shows that the following $\tilde g(z)$ does the job: 
\[ 
\begin{pmatrix}
\sum_{l=0}^{\omega-1} z^la^{(0)}_{1+nl,1}  &  z\sum_{l=0}^{\omega-1} z^la^{(-1)}_{n+nl,1}  & \hdots & z\sum_{l=0}^{\omega-1} z^la^{(-n+1)}_{2+nl,1}\\
\sum_{l=0}^{\omega-1} z^l a^{(0)}_{2+nl,1} & \sum_{l=0}^{\omega-1} z^la^{(-1)}_{1+nl,1} & \hdots & 
z\sum_{l=0}^{\omega-1} z^la^{(-n+1)}_{3+nl,1}\\
\vdots & \vdots & \ddots& &\\
\sum_{l=0}^{\omega-1} z^la^{(0)}_{n-1+nl,1} & \sum_{l=0}^{\omega-1} z^la^{(-1)}_{n-2+nl,1} & \hdots & z\sum_{l=0}^{\omega-1} z^la^{(-n+1)}_{n+nl,1}\\
\sum_{l=0}^{\omega-1} z^la^{(0)}_{n+nl,1} & \sum_{l=0}^{\omega-1} z^la^{(-1)}_{n-1+nl,1} & \hdots & 
\sum_{l=0}^{\omega-1} z^la^{(-n+1)}_{1+nl,1}
\end{pmatrix} 
\]
\end{proof}

\begin{cor}
For every $I\in \binom{[n]}{k}$ one has an equality
\[
\Psi {\rm Aut}_{\Delta_n}(U_{n\omega})p_I = {\mathfrak B} p_{w(I)}.
\]
\end{cor}
\begin{proof}
Thanks to Proposition \ref{AutB} the desired equality holds if we use the $GL_n$ version
of the Iwahori group instead of ${\mathfrak B}$. However, 
${\mathfrak B}p_{w(I)}$ coincides with the $GL_n$ Iwahori group orbit of the point 
$p_{w(I)}$.
\end{proof}

\begin{rem}
Corollary above claims that for any $A\in {\rm Aut}_{\Delta_n}(U_{n\omega})$ and any
$I\in\binom{[n]}{k}$ there exists an element $g(z)\in {\mathfrak B}$ such that $Ap_I=g(z)p_{w(I)}$. We note
that $g(z)$ does depend on $I$, but if we extend ${\mathfrak B}$ to the Iwahori subgroup
for $GL_n$ (as we do in Proposition \ref{AutB}), then $g(z)$ depends only on $A$.
\end{rem}

We summarize the discussion above as follows.
\begin{thm}
The quiver Grassmannian $X(k,n,\omega)$ is isomorphic to the union of Schubert varieties $X_{w(I)}\subset \eF_n$ for all $I\in\binom{[n]}{k}$.
\end{thm}

\begin{example}
We close this section with an example for $n=2$. The quiver Grassmannians we are interested
in are of the form $X(1,2,\omega)$.
The Weyl group $W_2$ is generated by $s_0$ and $s_1$; the elements of $W_2$ are of the form  $s_1s_0s_1\dots$ and
$s_0s_1s_0\dots$. In particular, for each $\omega>0$ there exist exactly two elements  $\sigma_1(\omega), \sigma_2(\omega)\in W_2$ of length $\omega$. 
The image of the quiver Grassmannian $X(1,2,\omega)$ inside the affine flag variety $\eF_2$ is equal to the union of the Schubert varieties
$X_{\sigma_1(\omega)}\cup X_{\sigma_2(\omega)}$. 
\end{example}

\subsection*{Acknowledgments} M.L acknowledges the PRIN2017 CUP E8419000480006, and the MatMod@TOV project of the "Department of Excellence" funded by MUR for the period 2023-2027 to the Department of Mathematics of the University of Rome ``Tor Vergata".

\appendix

\section{Flatness: an example}\label{sec:Flatness-example}
Let us consider the degeneration of $\bP^2={\rm Gr}(1,3)$ to the quiver Grassmannian $X(1,3)=X(1,3,1)$. It was claimed in \cite{FLP21} that this degeneration is not flat. More precisely, it was claimed in loc.cit. that the dimension of the degree $(1,1,1)$ component of the homogeneous coordinate ring of 
$X(1,3)$ is larger than $10$ (which is the dimension of the degree three component
of the homogeneous coordinate ring of ${\rm Gr}(1,3)$).
However, the argument contains a mistake (we note that according to \cite{G01} the degeneration
of any ${\rm Gr}(k,n)$ to $X(k,n)$ is flat).  Below we correct this mistake.  
 
The variety $X(1,3)$ sits inside $\bP^2\times \bP^2\times \bP^2$. The corresponding homogeneous coordinate ring is triply-graded. We  
denote the homogeneous
coordinates in the $i$-th $\bP^2$ by $x_i,y_i,z_i$ ($i=1,2,3$). In particular,
the degree $(1,1,1)$ homogeneous component of the coordinate ring sits
inside the polynomial ring in variables $x_\bullet,y_\bullet,z_\bullet$; it consists of polynomials
which are linear in each group of variables $x_i,y_i,z_i$ ($i=1,2,3$). 

Recall that $X(1,3)$ has three irreducible components; the open cells in these irreducible components are of 
the form
\[
\begin{pmatrix}
1 & 0 & 0\\ a_1 & 1 & 0\\ b_1 & a_1 & 1
\end{pmatrix},\quad
\begin{pmatrix}
0 & 1 & 0\\ 0 & a_2 & 1\\ 1 & b_2 & a_2
\end{pmatrix},\quad
\begin{pmatrix}
 0 & 0 & 1\\ 1 & 0 & a_3\\ a_3 & 1 & b_3
\end{pmatrix}.
\]
Here the $j$-th matrix corresponds to the $j$-th irreducible component ($j=1,2,3$);
for $i=1,2,3$ the $i$-th column of the $j$-th matrix produces the homogeneous coordinates 
$x_i,y_i,z_i$ of the points in the open cell of the $j$-th irreducible component
(we note that $a_\bullet$ and $b_\bullet$ are free parameters). 

The dimension of the degree three homogeneous component of the coordinate ring of the Grassmannian ${\rm Gr}(1,3)=\bP^2$ is equal to ten. Hence we have to show that the 
dimension of the degree $(1,1,1)$ homogeneous component of the coordinate ring of $X(1,3)$
is also equal to $10$. 
We claim that
the following ten monomials produce a basis of the degree $(1,1,1)$ component:
\begin{gather*}
y_1y_2z_3=(a_1,0,0), \quad y_1z_2y_3=(0,0,a_3),\quad  y_1z_2z_3=(a_1^2,0,b_3),\\
z_1y_2y_3=(0,a_2,0), \quad z_1y_2z_3=(b_1,a_2a_3,0), \quad z_1z_2y_3=(0,b_2,a_1a_3),\\
 z_1z_2z_3=(b_1a_2,b_2a_3,b_3a_1),  \\ x_1y_2z_3=(1,0,0), \quad
x_1y_3z_3=(0,1,0), \quad y_1z_2x_3=(0,0,1).
\end{gather*}
Here the right hand side of each equality computes the value of the corresponding monomial on the three open cells. The monomials above are clearly linearly independent. The values of all other
degree $(1,1,1)$ monomials are either zero or coincide with one of the above. For example, 
$x_1x_2z_3=0$ and $x_1z_2z_3=(a_1,0,0)=y_1y_2z_3$.

\end{document}